\documentclass[pdflatex,sn-mathphys]{sn-jnl}

\usepackage{amsmath}
\usepackage{amsfonts}
\usepackage{amsthm}
\usepackage{mathtools}
\usepackage{bm}
\DeclareSymbolFont{largesymbolsA}{U}{txexa}{m}{n}
\DeclareMathSymbol{\varprod}{\mathop}{largesymbolsA}{16}

\jyear{2022}%

\theoremstyle{thmstyleone}%
\theoremstyle{thmstyletwo}%

\theoremstyle{thmstylethree}%
\newtheorem{myalgorithm}{Algorithm}

\newcommand{\myR}{\mathbb{R}}
\newcommand{\myN}{\mathbb{N}}
\newcommand{\myRplus}{{\myR^+}}
\newcommand{\myNplus}{{\myN^+}}

\newcommand{\myPara}{\boldsymbol{\mu}}
\newcommand{\myParaTrue}{{\myPara^\star}}
\newcommand{\myParaOpt}{{\myPara_{\text{opt}}}}
\newcommand{\myParaDim}{{N_p}}
\newcommand{\myParaDomain}{\mathcal{P}}
\newcommand{\myParaSpace}{{\myR^\myParaDim}}

\newcommand{\myTimeDomain}{{I}}
\newcommand{\myTimeCoord}{{t}}

\newcommand{\myGeoDim}{{d}}
\newcommand{\myGeoSpace}{{\myR^\myGeoDim}}
\newcommand{\myGeoDomain}{{\Omega}}
\newcommand{\myGeoBorder}{{\partial \myGeoDomain}}
\newcommand{\myGeoCoord}{\mathbf{x}}

\newcommand{\myState}[3]{u(#1,#2;#3)}
\newcommand{\myStateRB}[3]{{u_\varepsilon(#1,#2;#3)}}
\newcommand{\myStateFO}[3]{{u_h(#1,#2;#3)}}
\newcommand{\mySolution}[1]{\myState{\myGeoCoord}{\myTimeCoord}{#1}}
\newcommand{\mySolutionRB}[1]{\myStateRB{\myGeoCoord}{\myTimeCoord}{#1}}
\newcommand{\mySolutionFO}[1]{\myStateFO{\myGeoCoord}{\myTimeCoord}{#1}}
\newcommand{\mySolutionSpace}{\mathcal{U}}

\newcommand{\solh}[1]{u_h({\myGeoCoord},{\myTimeCoord};{#1})} 
\newcommand{\solr}[1]{u_{\varepsilon}({\myGeoCoord},{\myTimeCoord};{#1})} 

\newcommand{\myObserver}{\mathcal{L}}
\newcommand{\myMeas}[1]{{\myObserver \mySolution{#1}}}
\newcommand{\myMeasRB}[1]{{\myObserver \solr{#1}}}
\newcommand{\myMeasDim}{{N_m}}
\newcommand{\myMeasSpace}{{\myR^\myMeasDim}}

\newcommand{\myForwardResponseMap}{\mathcal{G}}
\newcommand{\myForwardResponse}[1]{{\myForwardResponseMap(#1)}}
\newcommand{\myForwardResponseMapRB}{\mathcal{G}_\varepsilon}
\newcommand{\myForwardResponseRB}[1]{{\myForwardResponseMapRB(#1)}}

\newcommand{\myParticle}[2][j]{{\myPara_{#2}^{{}_{\left(#1\right)}}}}
\newcommand{\myParticleMean}[1]{{\overline{\myPara}_{#1}}}
\newcommand{\myParticleResponseMean}[1]{\overline{\myForwardResponseMap}_{#1}}
\newcommand{\myParticleResponseMeanRB}[1]{\overline{\myForwardResponseMap}_{{#1},\varepsilon}}

\newcommand{\myEnsemble}[1]{{\mathcal{E}_{#1}}}
\newcommand{\myEnsembleDim}{{J}}
\newcommand{\myEnsembleIndex}{{j}}
\newcommand{\myEnsembleExplicit}[1]{{\{\myParticle{#1}\}_{{}^{\myEnsembleIndex=1}}^{{}_\myEnsembleDim}}}
\newcommand{\myEnsembleDistr}[1]{{\Pi_{#1}}}

\newcommand{\myMin}[1]{\min_{#1}\,}
\newcommand{\myArgMin}[1]{{\text{arg}\myMin{#1}}}
\newcommand{\myInverse}[1]{{{#1}^{-1}}}
\newcommand{\myCostFunction}[2]{{\Phi (#1 \,\vert\, #2) }}
\newcommand{\myEuclidianNorm}[1]{{\| #1 \|_2}}
\newcommand{\myRange}[1]{\{1,\ldots,#1\}}

\newcommand{\myData}{\mathbf{y}}
\newcommand{\myDataDim}{\myMeasDim}

\newcommand{\myNoise}{\boldsymbol{\eta}}
\newcommand{\myNoiseCovariance}{\boldsymbol{\Sigma}}
\newcommand{\myNoiseCovarianceSpace}{{\myR^{\myDataDim \times \myDataDim}}}
\newcommand{\myNoiseNorm}[1]{{\| #1 \|^2_\myInverse{\myNoiseCovariance}}}

\newcommand{\myKalmanGain}[1]{{\mathbf{K}_{#1}}}
\newcommand{\myPCovariance}[1]{{\mathbf{P}_{#1}}}
\newcommand{\myQCovariance}[1]{{\mathbf{Q}_{#1}}}
\newcommand{\myNoisyData}[2][j]{{\myData_{#2}^{{}_{(#1)}}}}

\begin{document}

\title{A Reduced Basis Ensemble Kalman Method}


\author*[1]{\fnm{Francesco A. B.} \sur{Silva}}\email{f.a.b.silva@tue.nl}

\author[1]{\fnm{Cecilia} \sur{Pagliantini}}\email{c.pagliantini@tue.nl}

\author[2]{\fnm{Martin} \sur{Grepl}}\email{grepl@igpm.rwth-aachen.de}

\author[1]{\fnm{Karen} \sur{Veroy}}\email{k.p.veroy@tue.nl}

\affil*[1]{\orgdiv{Department of Mathematics and Computer Science}, \orgname{Eindhoven University of Technology}, \orgaddress{\city{Eindhoven}, \postcode{5600 MB}, \country{The Netherlands}}}

\affil[2]{\orgdiv{Institute of Geometry and Practical Mathematics}, \orgname{RWTH Aachen University}, \orgaddress{\city{Aachen}, \postcode{52056}, \country{Germany}}}


\abstract{In the process of reproducing the state dynamics of parameter dependent distributed systems, data from physical measurements can be incorporated into the mathematical model to reduce the parameter uncertainty and, consequently, improve the state prediction. Such a Data Assimilation process must deal with the data and model misfit arising from experimental noise as well as model inaccuracies and uncertainties. In this work, we focus on the ensemble Kalman method (EnKM), a particle-based iterative regularization method designed for \textit{a posteriori} analysis of time series. The method is gradient free and, like the ensemble Kalman filter (EnKF), relies on a sample of parameters or particle ensemble to identify the state that better reproduces the physical observations, while preserving the physics of the system as described by the best knowledge model.
We consider systems described by parameterized parabolic partial differential equations and employ model order reduction (MOR) techniques to generate surrogate models of different accuracy with uncertain parameters. Their use in combination with the EnKM involves the introduction of the model bias which constitutes a new source of systematic error. To mitigate its impact, an algorithm adjustment is proposed accounting for a prior estimation of the bias in the data. The resulting RB-EnKM is tested in different conditions, including different ensemble sizes and increasing levels of experimental noise. The results are compared to those obtained with the standard EnKF and with the unadjusted algorithm.}




\maketitle

\section{Introduction}\label{Introduction}

The problem of estimating model parameters of static and dynamical systems is encountered in many applications from earth sciences to engineering. In this work we focus on the parameter estimation of dynamical systems described by parameterized parabolic partial differential equations (pPDEs). Here, we assume that a limited and polluted knowledge of the solution is available at multiple time instances through noisy local measurements. 

For solving this kind of inverse problem, countless deterministic and stochastic methods have been proposed. Among them, a widely used technique is the so-called ensemble Kalman filter \cite{evensen2003}, a recursive filter employing a series of measurements to obtain improved estimates of the variables involved in the process. The idea of using the EnKF for reconstructing the parameters of dynamical systems traces back to \cite{Anderson2001, Lorentzen2001}, in which trivial artificial dynamics for the parameters was assumed to make the estimation possible. This was naturally accompanied by efforts for improving the performance of the method in terms of stability, by introducing covariance inflation \cite{Hamill2001, AndersonAnderson1999} and localization \cite{Hamill2001, Houtekamer2001}, and in terms of computational cost. Relevant to the latter have been the development of multi- level methods \cite{Hoel2016}, the use of model order reduction techniques \cite{pagani2016reduced}, and the introduction of further surrogate modeling techniques \cite{Popov2022}. The use of approximated models inevitably led to the study of the impact of model error on the EnKF \cite{Herschel2002, Mitchell2015}, alongside with other data assimilation methods \cite{Calvetti2018, Huttunen2007}.

Although ensemble Kalman methods were originally meant for sequential data assimilation, i.e., for real-time applications, they proved to be reliable also for asynchronous data assimilation \cite{Sakov2010}. The first paper proposing to adapt the EnKF to a retrospective data analysis was \cite{Skjervheim2007}. For analysis, the data are employed all at once at the end of an assimilation window, which is in common with a series of methods, e.g., variational methods \cite{Zhijin2001} such as 4D-VAR \cite{Thepaut1991} and other smoothers \cite{anderson2012optimal}. Compared to those approaches, the EnKF is particularly appealing since it does not require the computation of Fr\'{e}chet derivatives, a major complication for data assimilation algorithms.

In \cite{Iglesias_2013}, Iglesias et al. introduced what they called the {\it{ensemble Kalman method}}, an EnKF-based asynchronous data assimilation algorithm. Depending on the design of the algorithm, this method has connections to Bayesian data assimilation \cite{Schillings2018} and to maximum likelihood estimation \cite{Chen2012}. In particular, in the latter case, the method constitutes an ensemble-based implementation of so-called {\it{iterative regularization methods}} \cite{Kaltenbacher2008}. In the case of perfect models, the EnKM has already been analyzed in depth in \cite{Schillings2018, Evensen2018} and convergence and identifiability enhancements have been proposed in \cite{Wu2019, Iglesias2016}. Due to the iterative nature of the EnKM, dealing with high-dimensional parametric problems is often computationally challenging. In \cite{Gao2021} a multi-level strategy has been proposed to improve the computational performance of the method. 

In this work we propose an algorithm, called {\it{Reduced Basis Ensemble Kalman Method}} (RB-EnKM), that leverages the computational efficiency of surrogate models obtained with MOR techniques to solve asynchronous data assimilation problems via ensemble Kalman methods. The use of the EnKM allows us to avoid adjoint problems that are often difficult to reduce and intrinsically depend on the choice of measurement positions. Model order reduction, already employed in other data assimilation problems \cite{Gong2019, Nadal2015}, is used as a key tool for accelerating the method. However, the use of approximate models within the EnKM introduces a model error that could hinder the convergence of the method. In this work, we propose to deal with this error by including a prior estimation of the bias in the data.
Specifically, we incorporate empirical estimates of the mean and covariance of the bias in the Kalman gain. In some instances, those quantities can be computed at a negligible cost by employing the same training set used for the construction of the reduced model.

The paper is structured as follows: in Section \ref{sec:problem_formulation} we introduce the asynchronous data assimilation problem together with the standard ensemble Kalman method (Algorithm \ref{alg:enkm_algorithm}). Subsequently, in Section \ref{sec:surrogate_models}, we present an overview on reduced basis (RB) methods and describe how to use them in combination with the ensemble Kalman method to derive the RB-EnKM (Algorithm \ref{alg:rb_enkm_algorithm}). In Section \ref{sec:exp}, we test the new method on two numerical examples. In the first example, we estimate the diffusivity in a linear advection-dispersion problem in 2D (Section \ref{sec:linear_experiment}), while in the second, we estimate the hydraulic log-conductivity in a non-linear hydrological problem (Section \ref{sec:nonlinear_experiment}). In both cases, we compare the behavior of the full order and reduced order models in different conditions. Section \ref{sec:conclusions} provides conclusions and considerations on the proposed method and on its numerical performances.

\section{Problem Formulation}
\label{sec:problem_formulation}

Let $\mySolutionSpace$ be a given function space and let $\myParaDomain \subset \myParaSpace$, with $\myParaDim \in \myNplus$, be a set of model parameters. We consider the pPDE: for any parameter $\myPara \in \myParaDomain$, find $\myState{\,\cdot\,}{\,\cdot\;}{\myPara} \in \mathcal{U}$ such that $\partial_t \mySolution{\myPara} = \mathcal{F}_{\myPara} \mySolution{\myPara}$ for any $\myGeoCoord\in\myGeoDomain \subset \myGeoSpace$ and $t \in \myTimeDomain \coloneqq \left( 0, T \right] \subset \myRplus$. Here $\mathcal{F}_{\myPara}$ is a generic parameterized differential operator and $\partial_t$ is the first order partial time derivative. This pPDE provides the constraint to the inverse problem of estimating the unknown parameter $\myParaTrue \in \mathcal{P}$ from data or observations given by
\begin{align}
\begin{split}
    & \myData(\myParaTrue, \myNoise) = \myMeas{\myParaTrue} + \myNoise \\
    & s.t. \quad \partial_t \mySolution{\myParaTrue} = \mathcal{F}_{\myParaTrue} \mySolution{\myParaTrue}. \label{eq:pPDE}
\end{split}
\end{align}
Here, $\myObserver : \mySolutionSpace \rightarrow \myMeasSpace$, with $\myMeasDim \in \myNplus$, maps the space of the solutions to the space of the measurements, simulating the observation process, and $\myNoise$ is an unknown realization of a Gaussian random variable with zero mean and given covariance, $\myNoiseCovariance \in \myNoiseCovarianceSpace$. Note that both the observed data $\myData$ and additive noise $\myNoise$ are $\myDataDim$-dimensional vector-valued quantities and that $\myNoiseCovariance$ is a symmetric positive-definite matrix defining the inner product $\myNoiseNorm{\cdot} \coloneqq \myEuclidianNorm{ \myNoiseCovariance^{-1/2} \cdot}$ on $\myR^\myDataDim$, where $\myEuclidianNorm{\cdot}$ is the Euclidean norm. 

To solve this inverse problem, we must explicitly solve the pPDE~\eqref{eq:pPDE}. This is done using a suitable discretization, in space and time, of the differential operator $\mathcal{F}_{\myPara}$. To this end, we introduce an approximation space $ \mathcal{V}_h \subset \mathcal{U}$ so that the approximate problem reads:
\begin{equation}\label{eq:FOM_pPDE}
    \text{find} \;\;  u_h(\myPara)=u_h(\,\cdot\,,\,\cdot\,;\myPara) \in \mathcal{V}_h \quad \text{s.t.} \quad \partial_t \solh{\myPara} = \mathcal{F}^h_{\myPara}\solh{\myPara}.
\end{equation}
The discretization of the pPDE can be chosen according to the specific problem of interest. In all numerical examples proposed in this work, we employ a space-time Petrov--Galerkin discretization of \eqref{eq:pPDE} with piecewise polynomial trial and test spaces, as described in Section~\ref{sec:exp}, and we assume \eqref{eq:FOM_pPDE} to be sufficiently accurate such that we can take $\myData(\myParaTrue, \myNoise) = \mathcal{L}\mySolutionFO{\myParaTrue} + \myNoise$.

To characterize the observation of the solution, we introduce the \textit{forward response} map $\myForwardResponseMap : \myParaDomain \rightarrow \myMeasSpace$ defined as $\myForwardResponse{\myPara} \coloneqq \mathcal{L} \mySolutionFO{\myPara} $ for any solution of the pPDE \eqref{eq:FOM_pPDE}. Although the use of the map $\mathcal{G}$ results in a more compact notation, omitting its dependence on the solution of the pPDE conceals a key aspect of the method, i.e., the mapping from the parameter vector to the corresponding space-time pPDE solution. For this reason, and because it makes it harder to introduce the problem discretization, it will be used with caution.

\subsection{The Ensemble Kalman Method}
\label{sec:enkm}

The data assimilation problem presented above can be recast as a minimization problem for the cost functional, $\myCostFunction{\myPara}{\myData} \coloneqq \myNoiseNorm{ \myData (\myParaTrue, \myNoise) - \mathcal{L} \mySolutionFO{\myPara} }$, representing the misfit between the experimental data, $\myData (\myParaTrue, \myNoise)$, and the forward response. The optimal parameter estimate $\myParaOpt(\myData)$ is thus given by
\begin{align}
\begin{split}
    \label{eq:cost_functional}
    &\myParaOpt (\myData) = \myArgMin{\myPara \in \myParaDomain} \myCostFunction{\myPara}{\myData} \\
    &s.t. \quad \partial_t \mySolutionFO{\myPara} = \mathcal{F}^h_{\myPara} \mySolutionFO{\myPara}.
\end{split}
\end{align}
This is equivalent to a maximum likelihood estimation, given the {\it{likelihood}} function, $l(\myPara \,\vert\, \myData) = \exp\{-\frac{1}{2} \myCostFunction{\myPara}{\myData} \}$, associated with the probability density function of the data, $\myData \vert \myPara$, i.e., the probability of observing $\myData$ if $\myPara$ is the parametric state. The shape of the function follows from the probability density function of the Gaussian noise realization.

Among various methods proposed to solve this optimization problem, the EnKM relies on a sequence of parameter ensembles $\myEnsemble{n}$, with $n \in \myNplus$, to estimate the minimum of the cost functional. Each ensemble consists of a collection $\myEnsembleExplicit{n}$ of $\myEnsembleDim \in \myNplus$ parameter vectors $\myParticle{n}$, hereby named ensemble members or particles, whose interaction, guided by the experimental measurements, causes them to cluster around the solution of the problem as iterations proceed. 
At the beginning of each iteration, the solution of the pPDE and its observations are computed for each $\myEnsembleIndex \in \myRange{\myEnsembleDim}$. Subsequently, the ensemble is updated based on the empirical correlation among parameters and between parameters and measurements, as well as on the misfits between the experimental measurements $\myData (\myParaTrue, \myNoise)$ and the particle measurements $\mathcal{L} \mySolutionFO{\myParticle{n}}$. A single iteration, equivalent to the one in \cite{Iglesias_2013}, is formalized in the following pseudo algorithm:

\begin{myalgorithm} Iterative ensemble method for inverse problems.
\label{alg:enkm_algorithm}
\hfill \break

Let $\myEnsemble{0}$ be the initial ensemble with elements $\{\myParticle{0}\}_{{}^{j=1}}^{{}_J}$ sampled from a distribution $\myEnsembleDistr{0} ( \myPara )$.
For $n = 0,1,\ldots$

\hfill \break
\textit{(i)} \textbf{Prediction step}. Compute the measurements of the solution for each particle in the last updated ensemble:
\begin{align}
\begin{split}
    \myForwardResponse{\myParticle{n}} &= \mathcal{L} \mySolutionFO{\myParticle{n}} \\ s.t. \quad \partial_t \mySolutionFO{\myParticle{n}} &= \mathcal{F}^h_{\myParticle{n}} \mySolutionFO{\myParticle{n}} \quad \text{for all  } \myEnsembleIndex \in \myRange{\myEnsembleDim}. \label{eq:preliminary_step}
\end{split}
\end{align}

\hfill \break
\textit{(ii)} \textbf{Intermediate step}. From the last updated ensemble measurements and parameters, define the sample means and covariances:
\begin{align}
\mathbf{P}_n &= \frac{1}{\myEnsembleDim} \sum_{\myEnsembleIndex=1}^\myEnsembleDim \myForwardResponse{\myParticle{n}} \myForwardResponse{\myParticle{n}}^\top - \;\myParticleResponseMean{n} \myParticleResponseMean{n}^\top && \quad \text{with} \quad \myParticleResponseMean{n} = \frac{1}{\myEnsembleDim} \sum_{\myEnsembleIndex=1}^\myEnsembleDim \myForwardResponse{ \myParticle{n} }
\label{eq:prediction_step_P}\\
\mathbf{Q}_n &= \frac{1}{\myEnsembleDim} \sum_{\myEnsembleIndex=1}^\myEnsembleDim \myParticle{n} \myForwardResponse{\myParticle{n}}^\top - \;\myParticleMean{n} \myParticleResponseMean{n}^\top && \quad \text{with} \quad \myParticleMean{n} = \frac{1}{\myEnsembleDim} \sum_{\myEnsembleIndex=1}^\myEnsembleDim  \myParticle{n}. \label{eq:prediction_step_Q}
\end{align}

\hfill \break
\textit{(iii)} \textbf{Analysis step}. Update each particle in the ensemble:
\begin{align}
    \myParticle{n+1} = \myParticle{n} + \mathbf{Q}_n ( \mathbf{P}_n + \myNoiseCovariance )^{-1} (\myNoisyData{n} - \myForwardResponse{ \myParticle{n} } ) \label{eq:analysis_step} \quad \text{with} \quad \myNoisyData{n} \sim \mathcal{N} (\myData, \myNoiseCovariance).
\end{align}
\end{myalgorithm}

In the last step of the algorithm, the cross correlation matrices $\myPCovariance{n}$ and $\myQCovariance{n}$ are used to compute the Kalman gain $\myKalmanGain{n} \coloneqq \myQCovariance{n} \myInverse{ (\myPCovariance{n} + \myNoiseCovariance) }$. This modulates the extent of the correction: a low-gain corresponds to conservative behavior, i.e., small changes in the particle positions, while a high-gain involves a larger correction. Note that the experimental data are perturbed with artificial noise sampled from the same distribution assumed for the experimental noise $\myNoise$. This leads to an improved estimate over the unperturbed case.

A termination criterion for the algorithm is essential for the proper implementation of the method. The one presented in \cite{Iglesias_2013} is based on the discrepancy principle and consists in stopping the algorithm when the error between the experimental data and the measurements is comparable to the experimental noise, that is, when $\myNoiseNorm{\myData-\myForwardResponse{\myParticleMean{n}}} \leq \sigma \myNoiseNorm{\myNoise}$ for some $\sigma \geq 1$. An alternative approach is to set a threshold for the norm of the parameter update, i.e., to terminate the algorithm when $\myEuclidianNorm{\myParticleMean{n+1}-\myParticleMean{n}} \leq \tau \myEuclidianNorm{\myParticleMean{n+1}}$ for some $\tau \ll 1$. The latter criterion is more robust to model errors and is therefore used in our numerical experiments.

Equally important for the method is the choice of the distribution $\myEnsembleDistr{0}$ from which the initial ensemble (or first guess) $\myEnsemble{0}$ is sampled. In most of the cases, including those considered in our numerical experiments, the distribution $\Pi_0$ comes from an \textit{a priori} knowledge of the range of admissible parameters. In other scenarios, e.g., when the parameters live in an infinite-dimensional space, it may be necessary to define additional criteria on how to treat the parameter space. The initial ensemble plays a fundamental role in stabilizing the inverse problem. Indeed, it has been shown in \cite{Iglesias_2013} that all the ensembles generated by Algorithm \ref{alg:enkm_algorithm} are contained in the space spanned by the initial ensemble, that is
\begin{equation}
    \label{closed_space}
    \myEnsemble{n} \in \mathcal{A} \coloneqq \text{span}\,
    \myEnsembleExplicit{0} \quad \text{for all  } \, n \in \myNplus .
\end{equation}
Furthermore, in the mean-field limit, i.e., in the case of infinite particles, and assuming an affine relationship between parameters and synthetic measurements, the distribution $\myEnsembleDistr{0}$ plays the same role as the Tikhonov regularization in variational data assimilation, see \cite{Asch2016}. In particular, the stabilization term is given by $-\log_e \myEnsembleDistr{0}(\myPara)$.

The main sources of error of the EnKM are associated with the ensemble size and with the evaluation of $\mathcal{G}(\myParticle{n}) = \mathcal{L}\mySolutionFO{\myParticle{n}}$. Indeed, while the observation of the solution is accurate and computationally cheap to evaluate, due to the linearity of the operator, the accuracy in the computation of the pPDE solution intrinsically depends on the quality of the numerical discretization. 
High order numerical discretizations might require prohibitively large computational costs, especially if the pPDE \eqref{eq:FOM_pPDE} is solved for many values of the parameter and over long temporal intervals.

The other steps of Algorithm~\ref{alg:enkm_algorithm} involve the following operations: (i) the assembly of $\mathbf{P}_n$ and $\mathbf{Q}_n$ in \eqref{eq:prediction_step_P}-\eqref{eq:prediction_step_Q}, with computational complexity of order $\mathcal{O}(J N_m^2)$, and (ii) the inversion of the matrix $\mathbf{P}_n+\myNoiseCovariance$ in the \emph{analysis step} \eqref{eq:analysis_step} with complexity $\mathcal{O}(N_m^3)$. 
The solution of the pPDE \eqref{eq:FOM_pPDE}, $\text{for all  } \myEnsembleIndex \in \myRange{\myEnsembleDim}$, in the prediction step of Algorithm \ref{alg:enkm_algorithm} is thus the computational bottleneck of the EnKM algorithm.

\section{Surrogate Models}
\label{sec:surrogate_model_general}

\subsection{Reduced Basis Methods}
\label{sec:surrogate_models}

Given the need to solve the pPDE \eqref{eq:FOM_pPDE} for several instances of the parameter, the use of MOR techniques appears an ideal choice.
Model order reduction has allowed exceptional computational speed-ups
in settings that require repeated model evaluations, such as multi-query simulations. In MOR the high-dimensional problem is replaced with a surrogate model of reduced dimensionality that still possesses optimal or near-optimal approximation properties but that can be solved at a considerably reduced computational cost.
In this work, we focus on a particular class of MOR techniques, known as reduced basis methods \cite{Pru02}. 

The reduced basis method typically consists of two phases: an offline phase and an online phase. In the computationally expensive offline phase a low-dimensional approximation of the solution space, namely the \emph{reduced space}, is constructed and a surrogate model is derived via projection of the full order model onto the reduced space. Then, the resulting low-dimensional reduced model can be solved in the online phase for many instances of the parameter at a computational cost independent of the size of the full order model.

To be more precise, let $\mathcal{M} := \{u_h(\myPara) \in \mathcal{V}_h \subset\mySolutionSpace \;\vert\; \partial_t \solh{\myPara} = \mathcal{F}^h_{\myPara} \solh{\myPara} \; \mbox{for all}\, \myPara \in \mathcal{P} \}$ be the solution set which collects the solution of the discretized pPDE \eqref{eq:FOM_pPDE} under variation of the parameter $\myPara$.
The parametric problem \eqref{eq:FOM_pPDE} is said to be reducible if the solution set $\mathcal{M}$ can be well approximated by a low-dimensional linear subspace.
In this case, such a subspace is obtained as the span of a problem-dependent basis derived from a collection of full order solutions or snapshots, $\{ u_h(\myPara_s) \}_{s=1}^{S}$, with $S \in \myNplus$, at sampled values, $\{ \myPara_s \}_{s=1}^{S}$, of the parameter. The set $\myParaDomain_\text{TRAIN}:=\{ \myPara_s \}_{s=1}^{S}\subset\myParaDomain$ of training parameters is a sufficiently rich subset of the parameter space that can be obtained by drawing random samples from a uniform distribution in $\myParaDomain$ or with other sampling techniques, such as statistical methods and sparse grids, see \cite[Chapter 6]{quarteroni2015reduced} and references therein.
The extraction of the basis functions from the snapshots is usually performed using SVD-type algorithms such as the proper orthogonal decomposition (POD) \cite{berkooz1993proper} or greedy algorithms. For problems that depend on both time and parameters, the so-called POD-Greedy method \cite{GP05,HO08} combines a greedy algorithm in parameter space with the proper orthogonal decomposition in time at a given parameter.
In the numerical tests of this work, we rely on the Weak-POD-Greedy algorithm, which is the preferred method whenever a rigorous error bound can be derived, while the POD and the Strong-POD-Greedy are often used when a bound is unavailable, i.e., for most of non-linear problems. Note that, under the same choice of training parameters, the latter are more accurate but computationally less efficient.

Once an $N_\varepsilon$-dimensional reduced basis $\{ \psi_{i}\}_{i=1}^{N_\varepsilon}$ is obtained, the associated reduced space is given by $\mathcal{V}_{\varepsilon} = \text{span} \{ \psi_1,\ldots,\psi_{N_\varepsilon}\}\subset\mathcal{V}_h$, and the full model solution $u_h(\myPara)$, for a given $\myPara$, is approximated with a function  $u_{\varepsilon}(\myPara)$ in $\mathcal{V}_{\varepsilon}$,
\begin{equation*}
u_{\varepsilon}(\myPara) = \sum_{i=1}^{N_\varepsilon} u_{i}(\myPara)\,\psi_i\,,\qquad \myPara\in\myParaDomain,
\end{equation*}
where $(u_{1}(\myPara),\ldots,u_{N_\varepsilon}(\myPara))^\top\in\mathbb{R}^{N_\varepsilon}$ denotes the vector of the expansion coefficients in the reduced basis.
The reduced model thus reads:
\begin{equation}\label{eq:ROM_pPDE}
    \text{find} \;\;  u_{\varepsilon}(\myPara)=u_{\varepsilon}(\,\cdot\,,\,\cdot\,;\myPara) \in \mathcal{V}_{\varepsilon} \quad \text{s.t.} \quad \partial_t \solr{\myPara} = \mathcal{F}^{\varepsilon}_{\myPara}\solr{\myPara},
\end{equation}
where the operator $\mathcal{F}^{\varepsilon}_{\myPara}:\mathcal{V}_{\varepsilon}\rightarrow\mathbb{R}$ is obtained by projecting the full order operator $\mathcal{F}^{h}_{\myPara}:\mathcal{V}_h \rightarrow \mathbb{R}$ from \eqref{eq:FOM_pPDE} onto the reduced space $\mathcal{V}_{\varepsilon}$.

The computational gain derived from solving problem \eqref{eq:ROM_pPDE} instead of the full order model \eqref{eq:FOM_pPDE} hinges on the feasibility of a complete decoupling of the offline and online phases.
A computational complexity of the online phase independent of the size of the full order problem can be achieved under the assumption of linearity and parameter-separability of the operator $\mathcal{F}^{h}_{\myPara}$. To deal with general non-linear operators, hyper-reduction techniques are required. These include methods for approximating the high-dimensional non-linear term $\mathcal{F}^{h}_{\myPara}$ with an empirical affine decomposition, such as the EIM \cite{barrault2004}, and methods for reducing the cost of evaluating the non-linear term, such as linear program empirical quadrature \cite{yano2019lp} and empirical cubature \cite{hernandez2017}.

\subsection{A Reduced Basis Ensemble Kalman Method}
\label{sec:rb_enkm}

In this section, we discuss the implications of replacing the high-fidelity model in the prediction step of the EnKM by a surrogate model derived via model order reduction, as described in Section~\ref{sec:surrogate_models}. The use of MOR for particle-based methods is particularly desirable in multi-query contexts since it allows us to significantly reduce the computational cost of solving the inverse problem. However, the approximation introduced by the model order reduction inevitably produces (small) deviations of the reduced solution from the full order one. This constitutes a problem for data assimilation algorithms, as already documented and investigated in \cite{Calvetti2018} and in other works.
Indeed, the error in the solution results in discrepancies between approximated and exact measurements. Although we can expect the mismatch $\boldsymbol{\delta}_\varepsilon(\myPara) \coloneqq \mathcal{L}u_h(\myPara) - \mathcal{L}u_\varepsilon(\myPara)$ to decrease with the approximation error of $u_\varepsilon(\myPara)$, this bias will inevitably entail a distortion of the loss functional obtained by simple model substitution, i.e., 
\begin{equation}
    \widetilde{\Phi} (\myPara\vert\,\myData) \coloneqq \| \myData(\myParaTrue, \myNoise) - \mathcal{L}\mySolutionRB{\myPara}\|^2_{\myNoiseCovariance^{-1}}.
    \label{eq:biased_cost_function}
\end{equation}
Note that this cost function does not vanish in the parameter $\myPara$ we are trying to estimate, not even in noise free conditions. This systematic error, independent of the magnitude of the experimental noise, can be mitigated by modifying the cost function and consequently the EnKM. A modified algorithm, which we refer to as {\it{adjusted RB-EnKM}} is presented in the following sections. This algorithm is in contrast to what we refer to as the {\it{biased RB-EnKM}}, i.e., the algorithm obtained by the simple substitution of the full order model with the reduced order model in Algorithm~\ref{alg:enkm_algorithm}, as presented in \eqref{eq:biased_cost_function}.

The modification of the algorithm can proceed in two ways. One possibility is to rewrite the exact cost function in terms of the surrogate model and the measurement bias, namely substituting $\mathcal{L}\mySolutionFO{\myPara} = \mathcal{L}\mySolutionRB{\myPara} + \boldsymbol{\delta}_\varepsilon(\myPara)$ in the minimization problem \eqref{eq:cost_functional} to obtain
\begin{align}
\begin{split}
    \Phi_1 (\myPara \vert\,\myData) \coloneqq& \; \| \myData(\myParaTrue, \myNoise) -  \mathcal{L}\mySolutionRB{\myPara} - \boldsymbol{\delta}_\varepsilon(\myPara) \|^2_{\myNoiseCovariance^{-1}} \\
    =& \; \| \mathcal{L}\mySolutionFO{\myParaTrue} -  \mathcal{L}\mySolutionFO{\myPara}  + \myNoise \, \|^2_{\myNoiseCovariance^{-1}}
    \equiv \Phi (\myPara\vert\,\myData).
\end{split}
\end{align}
A second option is to correct the experimental data involved in the biased cost function \eqref{eq:biased_cost_function} so that, at least in noise free conditions, its minimum coincides with the minimum of the exact cost function. This means subtracting $\boldsymbol{\delta}_\varepsilon(\myParaTrue)$ instead of $\boldsymbol{\delta}_\varepsilon(\myPara)$, and results in the new cost function
\begin{align}
\begin{split}
    \Phi_2 (\myPara \vert\,\myData) \coloneqq& \; \| \myData(\myParaTrue, \myNoise) -  \mathcal{L}\mySolutionRB{\myPara} -  \boldsymbol{\delta}_\varepsilon(\myParaTrue) \|^2_{\myNoiseCovariance^{-1}} \\
    =& \; \| \mathcal{L}\mySolutionRB{\myParaTrue} - \mathcal{L}\mySolutionRB{\myPara} + \myNoise \, \|^2_{\myNoiseCovariance^{-1}}
    \not\equiv \Phi (\myPara\vert\,\myData).
\end{split}
\label{eq:bias_corrected_cost_function}
\end{align}
In noise free conditions, i.e., if $\myNoise = \mathbf{0}$, both cost functions vanish at the exact value $\myParaTrue$. Since the cost functions $\Phi_1$ and $\Phi_2$ are non-negative, the minimum attained in $\myParaTrue$ is necessarily also a global minimum.

In the following, we focus on the second approach. The reason is that the first approach requires the evaluation of the bias at all parameter values $\myPara \in \mathcal{P}$ which is too expensive to perform. Furthermore, at the algorithmic level, the substitution of the true model with the sum of the surrogate model and its bias would significantly change the computation of $\mathbf{P}_n$ and $\mathbf{Q}_n$, and thus the algorithm structure. By contrast, the second approach is based on the assumption that the true model is incorrect and on the subsequent correction of the experimental data. This implies that $\boldsymbol{\delta}_\varepsilon(\myParaTrue)$ is the only bias involved and it requires just a single full order evaluation. However, since the argument $\myParaTrue$ is unknown, this is clearly not possible, and we must instead exploit the prior epistemic uncertainty on $\myParaTrue$, encoded in $\Pi_0(\myPara)$, to modify the cost function.

If $\myParaTrue$ is treated as a random variable with probability measure $\Pi_0$, then the data bias $\boldsymbol{\delta}^\star_\varepsilon = \boldsymbol{\delta}_\varepsilon(\myParaTrue)$ is in turn a random variable with probability measure $\Pi_0 \circ \boldsymbol{\delta}_\varepsilon^{-1}$. The moments of this distribution, henceforth denoted by $\overline{\boldsymbol{\delta}}_\varepsilon$ and $\boldsymbol{\Gamma}_\varepsilon$, can be empirically estimated via pointwise evaluations of the bias  without further assumptions on the nature of the distribution itself. However, the assumption of Gaussianity, although improperly implying the linearity of $\boldsymbol{\delta}_\varepsilon : \mathcal{P} \rightarrow \mathbb{R}^{N_m}$, is consistent with the other assumptions of Gaussianity and linearity required for the derivation of the EnKF \cite{evensen2003}. Furthermore, it allows us to obtain closed-form results, as shown in the next paragraphs.

In view of the fact that $\boldsymbol{\delta}^\star_\varepsilon$ is considered as a random variable, we change Equation \eqref{eq:bias_corrected_cost_function} to make explicit the dependence of the cost function $\Phi_2$ on $\boldsymbol{\delta}^\star_\varepsilon$; thereby
\begin{align}
    \Phi_\varepsilon (\myPara \, \vert \, \myData, \boldsymbol{\delta}^\star_\varepsilon) \coloneqq \; \|  \myData(\myParaTrue, \myNoise)  - \mathcal{L}\mySolutionRB{\myPara} -   \boldsymbol{\delta}^\star_\varepsilon\|^2_{\myNoiseCovariance^{-1}}. \label{eq:two_random_variable_loss_function}
\end{align}
In order to make the estimate of $\myPara$ dependent only on the experimental data, we must remove the conditioning on $\boldsymbol{\delta}^\star_\varepsilon$, i.e., marginalize out the random variable. The easiest way to do this is by employing Bayesian statistics and particularly recovering the same marginal distribution $\myData \vert \myPara$ mentioned at the beginning of Section~\ref{sec:enkm}. To this end, we consider the likelihood function $ l (\myPara \, \vert \, \myData, \boldsymbol{\delta}^\star_\varepsilon) \coloneqq \exp \{  - \frac{1}{2} \Phi_\varepsilon (\myPara \, \vert \,
\myData , \boldsymbol{\delta}^\star_\varepsilon) \}$, proportional to the density of $(\myData \,\vert\, \myPara, \boldsymbol{\delta}^\star_\varepsilon) \sim \mathcal{N}( \boldsymbol{\delta}^\star_\varepsilon + \mathcal{L}\mySolutionRB{\myPara}, \myNoiseCovariance)$. Employing \cite[Lemma 1.A]{sarkka_2013}, concerning the mean and covariance of the joint distribution of Gaussian variables, it can be easily proven that, if $\boldsymbol{\delta}^\star_\varepsilon \sim \mathcal{N}(\overline{\boldsymbol{\delta}}_\varepsilon, \Gamma_\varepsilon)$, then $\myData \,\vert\, \myPara \sim \mathcal{N}( \overline{\boldsymbol{\delta}}_\varepsilon + \mathcal{L}\mySolutionRB{\myPara}, \myNoiseCovariance + \Gamma_\varepsilon)$ and consequently we derive the marginalized cost functional
\begin{align}
    \Phi_\varepsilon (\myPara \, \vert \, \myData) \coloneqq \; \|  \myData - \mathcal{L}\mySolutionRB{\myPara} - \overline{\boldsymbol{\delta}}_\varepsilon \|^2_{(\myNoiseCovariance+\Gamma_\varepsilon)^{-1}}.
    \label{eq:marginalised_loss_function}
\end{align}

Hence, by analogy with Section~\ref{sec:enkm}, we can adapt the EnKM to optimize the new cost function under the surrogate model constraint \eqref{eq:ROM_pPDE}. The resulting adjusted RB-EnKM is summarized in Algorithm \ref{alg:rb_enkm_algorithm}. Unlike the reference EnKM, we distinguish between an offline and an online phase. In the offline phase, the training set of full order solutions is generated and used both to construct the surrogate model and to estimate the moments of $\boldsymbol{\delta}^\star_\varepsilon$. In the online phase, the actual optimization is performed.

\vspace{0.2cm} 
\begin{myalgorithm} Iterative ensemble method with reduced basis surrogate models and accounting for the associated measurements bias.
\label{alg:rb_enkm_algorithm}

\hfill \break
\textbf{Offline}:
Let $\mathcal{P}_{\text{TRAIN}}$ be a set of $S$ parameters $\{\myPara_{0}^{{}_{(s)}}\}_{{}^{s=1}}^{{}_S}$ sampled from the distribution $\myEnsembleDistr{0} ( \myPara )$. Based on the associated full order and reduced solutions $u_h(\myPara^{{}_{(s)}})$ and $u_\varepsilon(\myPara^{{}_{(s)}})$, respectively, we define the training biases
\begin{align}
    \boldsymbol{\delta}_{\varepsilon}(\myPara^{{}_{(s)}}) = \mathcal{L} u_h({\myPara}^{{}_{(s)}}) - \mathcal{L} u_{\varepsilon}({\myPara}^{{}_{(s)}}) \label{eq:meas_bias}
\end{align}
and the associated empirical moments
\vspace{-0.1cm}
\begin{align}
\boldsymbol{\Gamma}_{\varepsilon} = \frac{1}{S} \sum_{s=1}^S \boldsymbol{\delta}_{\varepsilon}(\myPara^{{}_{(s)}}) \boldsymbol{\delta}_{\varepsilon}(\myPara^{{}_{(s)}})^\top - \;\overline{\boldsymbol{\delta}}_\varepsilon \overline{\boldsymbol{\delta}}_\varepsilon^\top \quad \text{with} \quad \;\;\overline{\boldsymbol{\delta}}_\varepsilon = \frac{1}{S} \sum_{s=1}^S  \boldsymbol{\delta}_{\varepsilon}(\myPara^{{}_{(s)}}).
\end{align}

\hfill \break
\textbf{Online}:
Let $\myEnsemble{0}$ be the initial ensemble with elements $\{\myParticle{0}\}_{{}^{j=1}}^{{}_J}$ sampled from the distribution $\myEnsembleDistr{0} ( \myPara )$.
For $n = 0,1, \ldots$

\hfill \break
\textit{(i)} \textbf{Prediction step}. Compute the biased measurements of the approximated solution for each particle in the last updated ensemble:
\begin{align}
\begin{split}
    \myForwardResponseRB{\myParticle{n}} &= \myMeasRB{\myParticle{n}} \\ s.t. \quad \partial_t \mySolutionRB{\myParticle{n}} &= \mathcal{F}^\varepsilon_{\myParticle{n}} \mySolutionRB{\myParticle{n}} \quad \text{for all  } \myEnsembleIndex \in \myRange{\myEnsembleDim}. \label{eq:preliminary_step_rb}
\end{split}
\end{align}
\hfill \break
\textit{(ii)} \textbf{Intermediate step}. From the last updated ensemble measurements and parameters, define the sample means and covariances:
\begin{align}
\mathbf{P}_{n,\varepsilon} &= \frac{1}{\myEnsembleDim} \sum_{\myEnsembleIndex=1}^\myEnsembleDim \myForwardResponseRB{\myParticle{n}}\, \myForwardResponseRB{\myParticle{n}}^\top - \;\myParticleResponseMeanRB{n} \, \myParticleResponseMeanRB{n}^\top && \; \text{with} \;\;\; \myParticleResponseMeanRB{n} = \frac{1}{\myEnsembleDim} \sum_{\myEnsembleIndex=1}^\myEnsembleDim \myForwardResponseRB{ \myParticle{n} },
\label{eq:prediction_step_P_rb}\\
\mathbf{Q}_{n,\varepsilon} &= \frac{1}{\myEnsembleDim} \sum_{\myEnsembleIndex=1}^\myEnsembleDim \myParticle{n} \myForwardResponseRB{\myParticle{n}}^\top - \;\myParticleMean{n} \myParticleResponseMeanRB{n}^\top && \; \text{with} \;\;\;\;\;\myParticleMean{n} = \frac{1}{\myEnsembleDim} \sum_{\myEnsembleIndex=1}^\myEnsembleDim  \myParticle{n}. 
\label{eq:prediction_step_Q_rb}
\end{align}

\hfill \break
\textit{(iii)} \textbf{Analysis step}. Update each particle in the ensemble:
\begin{align}
\begin{split}
    \myParticle{n+1} &= \myParticle{n} + \mathbf{Q}_{n,\varepsilon} \left( \mathbf{P}_{n,\varepsilon} + \boldsymbol{\Gamma}_\varepsilon + \boldsymbol{\Sigma} \right)^{-1} (\myNoisyData{n} - \myForwardResponseRB{ \myParticle{n} } ) \label{eq:analysis_step_fb},\\\myNoisyData{n} &\sim \mathcal{N} (\myData - \overline{\boldsymbol{\delta}}_\varepsilon, \myNoiseCovariance + \boldsymbol{\Gamma}_\varepsilon).
\end{split}
\end{align}
\end{myalgorithm}
The prior probability $\Pi_0 (\myPara)$ used for the estimation of the moments of $\boldsymbol{\delta}^\star_\varepsilon$ could be substituted at every iteration by an updated probability measure of $\myPara$. However, the computation of the updated probability measure might compromise the computational gain obtained with the use of reduced models. One possibility to address this shortcoming is to use a Gaussian process regression of the initial ensemble biases to estimate the moments of $\boldsymbol{\delta}^\star_\varepsilon$ with respect to the new probability measure of $\myPara$. The development and study of this strategy together with its effect on the accuracy and performances of the RB-EnKM will be investigated in future studies.

\section{Numerical Experiments}
\label{sec:exp}

In the following section, we consider two data assimilation problems for the estimation of model parameters in parabolic partial differential equations. The first problem involves a linear advection dispersion problem with unknown P\'eclet number. The corresponding model is linear in the observed state $c(\mu)$, but it is non-linear in the parameter to estimate map. The second problem concerns the transport of a contaminant in an unconfined aquifer with unknown hydraulic conductivity. It involves two coupled partial differential equations: a stationary non-linear equation which describes the pressure field induced by an external pumping force and a time-dependent linear equation describing the advection-dispersion of the contaminant in a medium whose properties depend non-linearly on the pressure field. 

Both models describe 2D systems, and each exhibits ideal characteristics to test the proposed algorithms. The first, while leading to a non-linear inverse problem, is sufficiently simple to allow for a comparison between the adjusted and biased RB-EnKM and the reference full order EnKM. Moreover, its affine dependence on the parameter enables the use of error bounds for the efficient construction of the reduced space. The second problem, which is non-linear and non-affine in the six-dimensional parameter vector, is complex enough to serve as a non-trivial challenge for the proposed RB-EnKM algorithm, while the reference EnKM cannot even be tested due to the computational cost. From an \textit{a priori} estimate, performing full order tests with the same statistical relevance as the reduced basis ones would have taken up to 20 days on our machine.

The two problems are presented in Section \ref{sec:linear_experiment} and Section \ref{sec:nonlinear_experiment}.
We first introduce the pPDE, then present the full order discretization followed by the reduced basis approximation. The measurement operator is then introduced, and a first analysis of the inversion method is carried out. Finally, we study the impact of the ensemble size, of the experimental noise magnitude, and of the error of the reduced model on the reconstruction error of the EnKM. All the computations are performed using Python on a computer with 2.20 GHz Intel Core i7-8750H processor and 32 GB of RAM.

\subsection{Taylor--Green Vortex Problem}
\label{sec:linear_experiment}

\begin{figure}[!b]
    \begin{center}
        \includegraphics[scale=0.78]{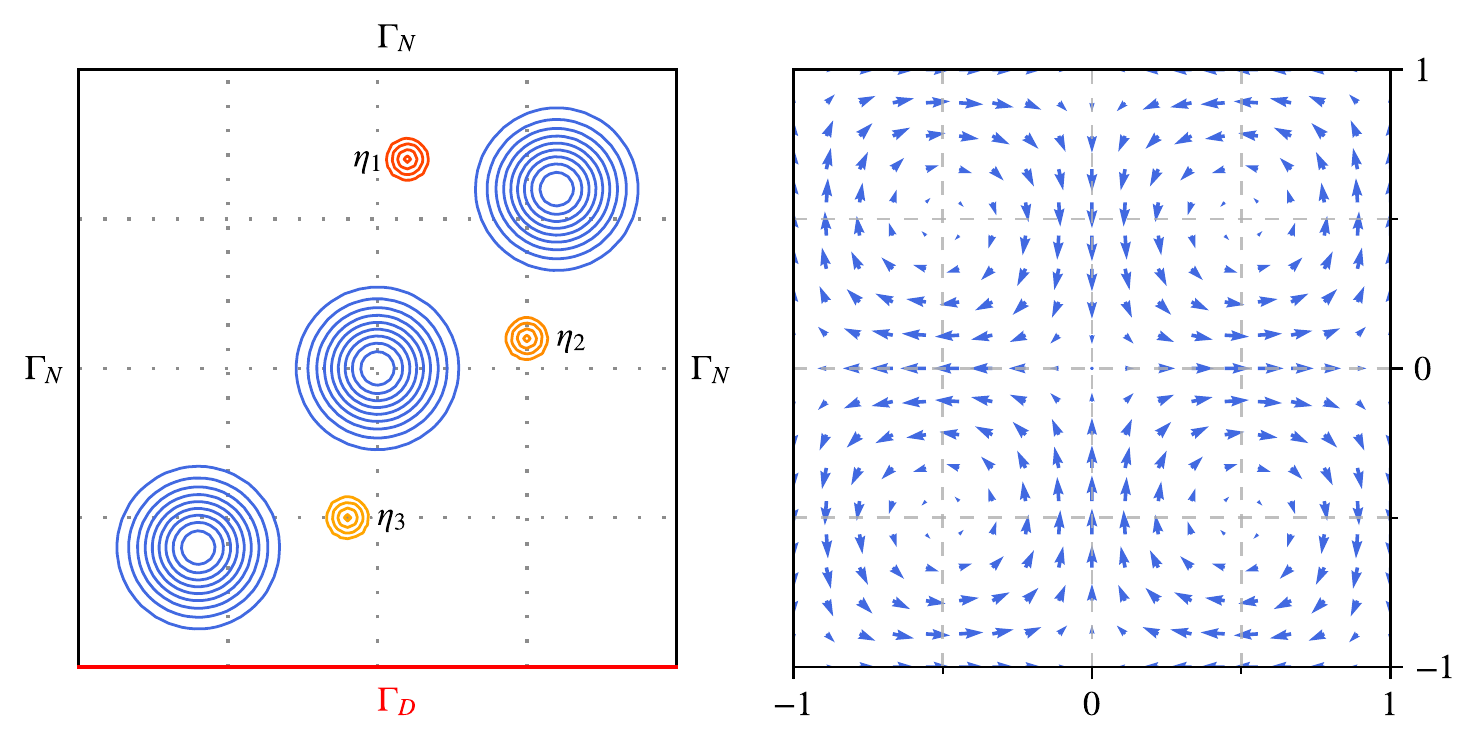}
    \end{center}
    \caption{Spatial domain of the Taylor--Green problem. On the left: the initial condition $c_0$ in blue, the sensor shape functions $\eta_i$, and the Neumann and Dirichlet boundaries, $\Gamma_N$, $\Gamma_D$. On the right, the velocity field, $\boldsymbol{\beta}$, with four Taylor--Green vortices.}
    \label{fig:domain_Taylor-Green}
\end{figure}

Let us consider the dispersion of a contaminant modeled by the 2D advection-diffusion equation with a Taylor--Green vortex velocity field \cite{karcher2018reduced}. 
We introduce the spatial domain $\myGeoDomain = (-1, 1)^2$ with Dirichlet boundary ${ \Gamma_D \coloneqq (-1,1) \times \{ -1 \} }$ and Neumann boundary ${ \Gamma_N \coloneqq \myGeoBorder \setminus \Gamma_D }$. 
We consider the problem of estimating the inverse of the P\'eclet number $\mu = 1/{\rm Pe}$ in the interval $\mathcal{P}:=[1/50, 1/10]$. The governing pPDE is given by: find $ c(\,\cdot\,, \,\cdot\,; \mu ) : \myGeoDomain \times \left(0,T\right] \rightarrow \mathbb{R} $ such that
\begin{equation}
\label{eq:linear_advection_diffusion}
\left\{
\begin{alignedat}{3}
    &\partial_t c - \mu \Delta c + \boldsymbol{\beta} \cdot \nabla c = 0, \qquad
        && \mbox{in} \; \myGeoDomain   & \times \, (0,T]\,,\\
    & \nabla c(\mathbf{x},t;\mu)\cdot \mathbf{n} = 0, 
        && \mbox{on} \; \Gamma_N & \times \, (0,T]\,,\\
    & c(\mathbf{x},t;\mu)=0, 
        && \mbox{on} \; \Gamma_D & \times \, (0,T]\,,\\
    & c(\mathbf{x},0;\mu)=c_0(\mathbf{x};\mu), 
        && \mbox{in} \; \myGeoDomain\,.&
\end{alignedat}
\right.    
\end{equation}
Here, the  velocity field $\boldsymbol{\beta} \coloneqq {(\sin(\pi x_1) \cos(\pi x_2)}$,$ {-\cos(\pi x_1) \sin(\pi x_2))^\top}$, $\mathbf{x} = (x_1, x_2)$, is a solenoidal field, and the initial condition $c_0(\,\cdot\,;\mu): \myGeoDomain \rightarrow \myR$ is given by the sum of three Wendland functions $\psi_{2,1}$ \cite{wendland1995piecewise} of radius $0.4$ and centers located at $(-0.6, -0.6)$, $(0, 0)$, and $(0.6, 0.6)$. The velocity field and the initial condition are shown in Figure \ref{fig:domain_Taylor-Green}. For the data assimilation problem we consider the time interval $\mathcal{I} \coloneqq \left( 0, 2.5 \right]$. 

The full order model is obtained by a nodal finite element discretization of Equation \eqref{eq:linear_advection_diffusion} using piecewise continuous polynomial functions, $\zeta_i: \myGeoDomain \rightarrow \myR$, $i = 1, \ldots, N_h$, of degree $2$ over a uniform Cartesian grid of width $h=0.04$, for a total of $N_h = 10,100$ degrees of freedom.
The resulting system of ordinary differential equations is integrated over time using a Crank--Nicolson scheme with uniform time step $\Delta t = 0.01$.
This is equivalent to performing a Petrov--Galerkin projection of Equation \eqref{eq:linear_advection_diffusion} with trial and test spaces defined as follows: we consider the partition of the temporal interval $\mathcal{I}$ into the union of equispaced subintervals, $\mathcal{I}_n \coloneqq \left(t_{n-1}, t_n \right]$, of length $\Delta t$ with $n=1,\ldots, N_t$ and $N_t \coloneqq 2.5 / \Delta t$.
Let $\omega_n :\mathcal{I} \rightarrow \myR$ be a piecewise constant function with support in $\mathcal{I}_n$, and let $\upsilon_n :\mathcal{I} \rightarrow \myR$ be a hat function with support in $\mathcal{I}_{n} \cup \mathcal{I}_{n+1}$. We define the trial space $\mathcal{V}_h \coloneqq \text{span} \{ \upsilon_n \cdot \zeta_i \}_{ {}^{i,n=1} }^{ {}_{N_h, N_t} }$ and the test space $\mathcal{W}_h \coloneqq \text{span} \{ \omega_n \cdot \zeta_i \}_{ {}^{i,n=1} }^{ {}_{N_h, N_t} }$, respectively.

\begin{figure}[!t]
    \begin{center}
        \includegraphics[scale=0.75]{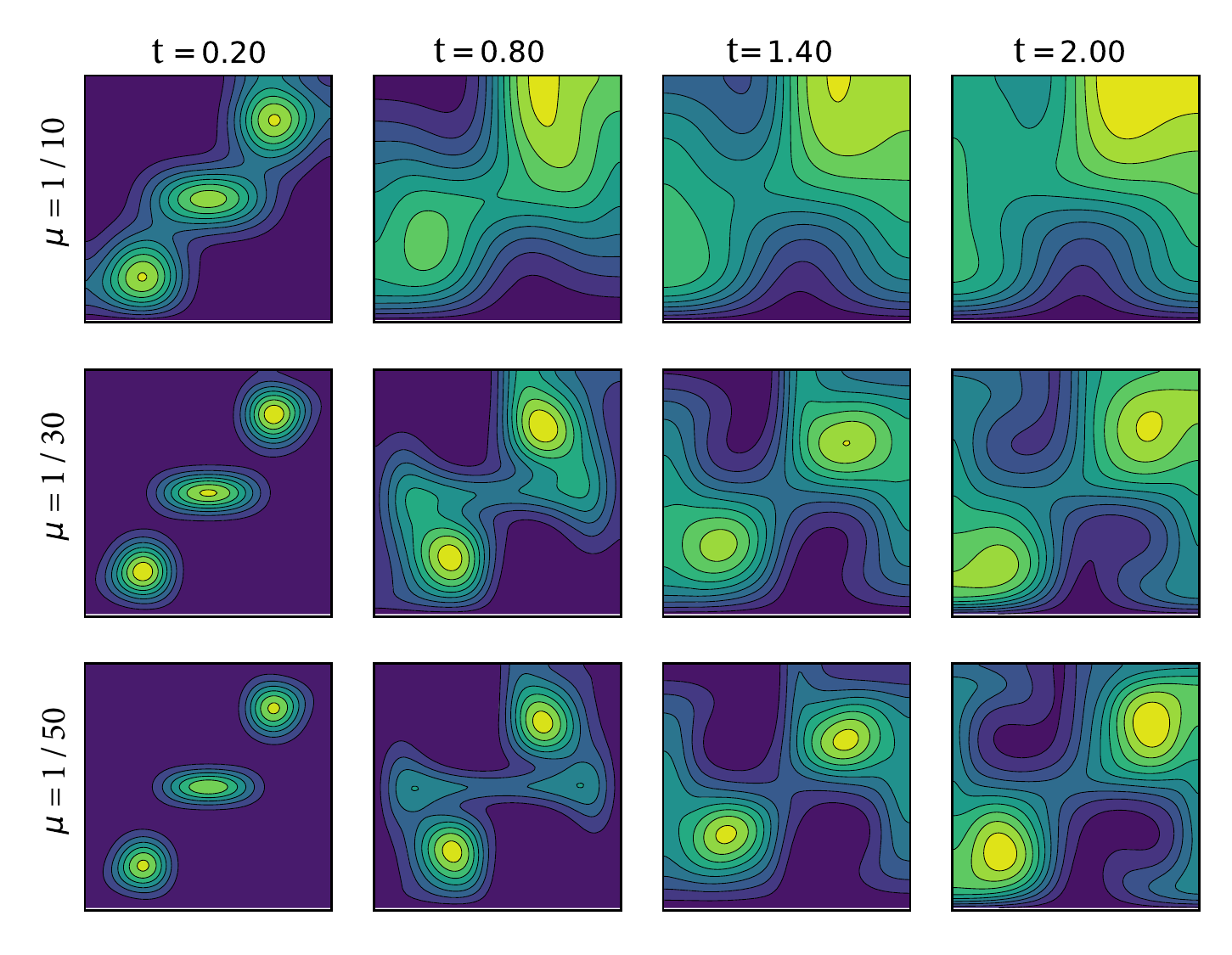}
    \end{center}
    \caption{Solution of the advection-diffusion equation for three increasing values of ${\rm Pe}$ at four time instances $t$. Snapshots normalized to unitary $L^\infty(\myGeoDomain)$ norm.} 
    \label{fig:solution_Taylor-Green}
\end{figure}

To solve the spatial problems arising at each time step, we use the sparse $\texttt{splu}$ function implemented in the $\texttt{scipy.sparse.linalg}\footnote{https://docs.scipy.org/doc/scipy/reference/sparse.linalg.html}$ package.
The computational time to obtain a single full order solution is on average $0.56$s. Snapshots of the solution at times $t \in \{0.2, 0.8, 1.4, 2.0\}$ and for the three parameter values $\mu \in \{1/10, 1/30, 1/50 \}$ are shown in Figure \ref{fig:solution_Taylor-Green}.

The high-fidelity model is used in combination with the time-gradient error bound $\Delta^{pr}_\text{R}(\mu)$ introduced in \cite{Nicole2022} to implement a Weak-POD-Greedy algorithm for the selection of the reduced basis functions.
To this end, we consider the training set $\Xi_{\text{TRAIN}}^\mu$ with parameters $\mu^{{}_{(s)}}= 1/(9.5 + 0.5 s)$ for all $s\in \mathbb{N}\cap [1,S]$ of size $S=81$. We prescribe a target accuracy of $10^{-2}$ for the maximum time-gradient relative error bound and we obtain an RB space of size $42$.
We can construct surrogate models of different accuracy by selecting  $N_\varepsilon \in \myN$ basis functions $\psi_i : \myGeoDomain \rightarrow \myR$, for $i=1,\ldots,N_{\varepsilon}$, out of these 42. Each choice corresponds to a relative error for the model given by
\begin{align}
    \label{eq:Taylor_Green_accuracy_definition}
    \varepsilon_c \coloneqq \sup_{\myPara \in \mathcal{D}} \frac{\| c_h(\myPara) - c_\varepsilon(\myPara) \|_{L^2(\mathcal{I},H^{1}(\myGeoDomain))} }{ \| c_h(\myPara) \|_{L^2(\mathcal{I},H^{1}(\myGeoDomain))} }.
\end{align}

Once the reduced basis has been computed, we construct a reduced model via a Petrov--Galerkin projection of Equation \eqref{eq:linear_advection_diffusion}
in the same way as we did for the full order model. For this purpose, we define the trial space $\mathcal{V}_\varepsilon \coloneqq \text{span} \{ \upsilon_n \cdot \psi_i \}_{ {}^{i,n=1} }^{ {}_{N_\varepsilon, N_t} }$ and the test space  $\mathcal{W}_\varepsilon \coloneqq \text{span} \{ \omega_n \cdot \psi_i \}_{ {}^{i,n=1} }^{ {}_{N_\varepsilon, N_t} }$.

We then look for a reduced solution of the form
\begin{equation}
    c_{\varepsilon}(\mathbf{x},t;\mu)=
    \sum_{i=1}^{N_{\varepsilon}}\sum_{n=1}^{N_t}
    c_{n,i}(\mu) \upsilon_n(t)\psi_i(\mathbf{x})\quad \mbox{for}\:t\in\mathcal{I},\;\mathbf{x}\in \Omega,
\end{equation}
where the expansion coefficients $c_{0,i}$, for $i=1,\ldots,N_\varepsilon$, result from the projection of the initial condition onto $\mathcal{V}_\varepsilon$, while the remaining coefficients $c_{n,i}$, with $i=1,\ldots,N_{\varepsilon}$ and $n=1,\ldots,N_t$, satisfy the equation
\begin{align}
    \label{eq:alg_concentration}
    \sum_{j=1}^{ N_\varepsilon} \left( \mathbf{M}_{ij} + \frac{\Delta t}{2} ( \mathbf{A}_{ij} + \mu \mathbf{K}_{ij}) \right) c_{n,i} = \sum_{j=1}^{ N_\varepsilon} \left( \mathbf{M}_{ij} - \frac{\Delta t}{2} ( \mathbf{A}_{ij} + \mu \mathbf{K}_{ij}) \right) c_{n-1,i}.
\end{align}
Here the matrices $\mathbf{M}, \mathbf{K}, \mathbf{A} \in \myR^{\scriptscriptstyle  N_\varepsilon \times N_\varepsilon}$ denote the mass, stiffness, and advection matrix, respectively, and are given by
\begin{equation}\label{eq:matrices}
    \mathbf{M}_{ij} \coloneqq \int_\Omega \psi_j \psi_i \, d\Omega,\;\;\, \mathbf{K}_{ij} \coloneqq \int_\Omega \nabla \psi_j \cdot \nabla \psi_i \, d\Omega,\;\;\, \mathbf{A}_{ij} \coloneqq \int_\Omega (\boldsymbol{\beta} \cdot \nabla \psi_j) \psi_i \, d\Omega.
\end{equation}

The solution of the system of equations \eqref{eq:alg_concentration}, equivalent to a Crank--Nicolson scheme, can be obtained iteratively solving $N_t$ linear systems of size $N_\varepsilon$ for an online complexity $\mathcal{O}(N_\varepsilon^3 + N_t N_\varepsilon^2)$. Employing all basis functions, the computational time for a reduced basis solution (online cost) is on average $5.4$ms, significantly less than the approximately $0.56$s required for a full order solution. The acceleration achieved is over $100$, which justifies the $47$s necessary for the construction of the RB model (offline cost), considering that the online phase requires computing up to 150 reduced basis solutions per iteration. Let us remark that such a cheap training phase is due to the low-dimensionality of the parameters space $\mathcal{P}$ and the availability of a tight error bound for this class of linear problems. 

\begin{figure}[!t]
    \begin{center}
        \includegraphics[scale=0.72]{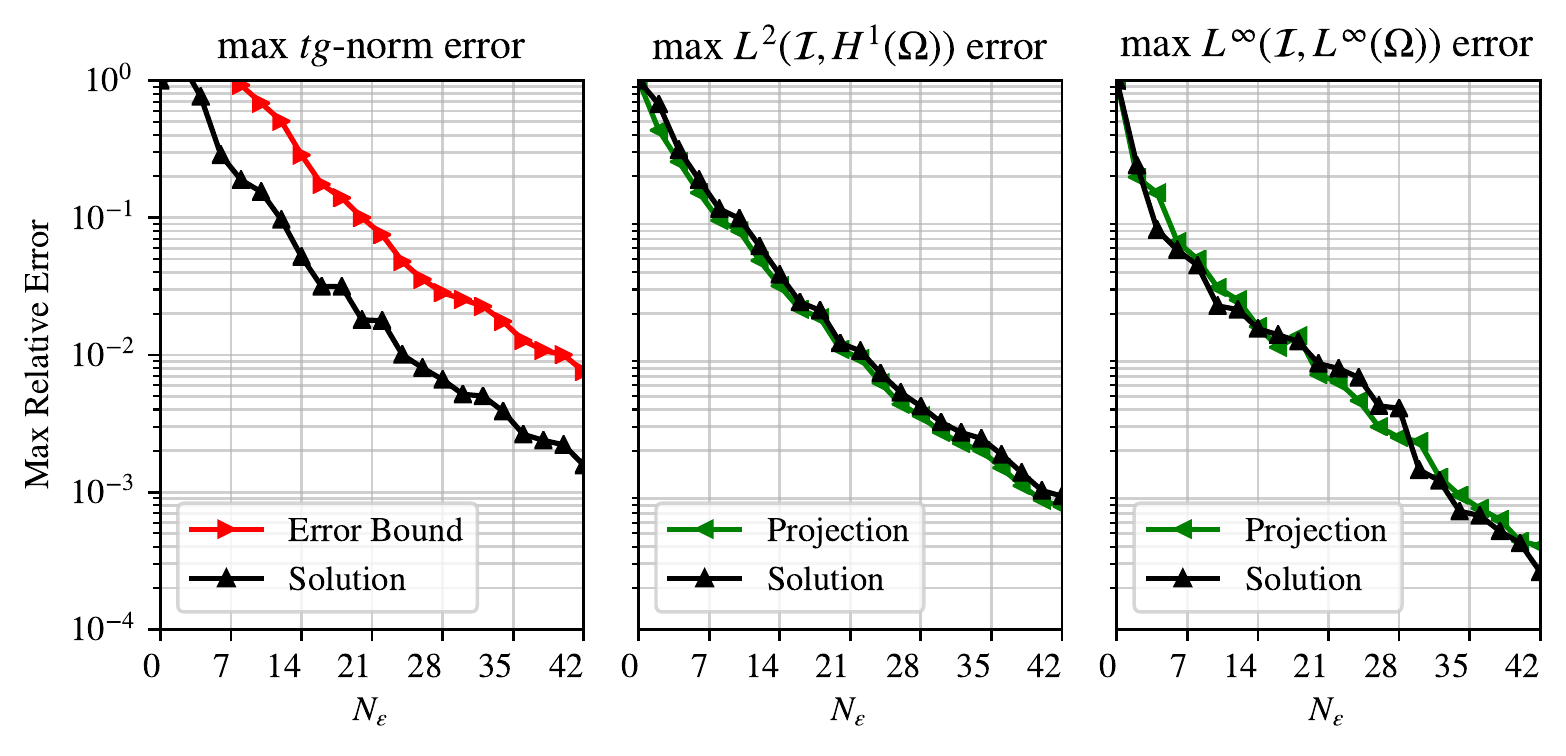}
    \end{center}
    \caption{Left: maximum relative time-gradient error and error bound of the advection-diffusion solution vs. $N_\varepsilon$. 
    Center: maximum $L^2(\mathcal{I}, H^1(\Omega))$ relative error of the projection and of the solution vs. $N_\varepsilon$. 
    Right: maximum $L^\infty(\mathcal{I}, L^\infty(\Omega))$ relative error of the projection and of the solution vs. $N_\varepsilon$. Projections based on the $L^2(\Omega)$ inner product of the gradients.}
    \label{fig:Taylor-Green_surrogate_error}
\end{figure}

Note that both the online computational cost and the accuracy of the solution depend on $\Delta t$ and on $N_\varepsilon$. The first is kept fixed, $\Delta t = 0.01$, while the latter varies in some of the experiments. In order to keep track of the error associated with different choices of $N_\varepsilon$, we proceed with the characterization of the error between the surrogate model solution $c_\varepsilon(\mu)$ and the full model solution $c_h (\mu)$ for different values of $N_\varepsilon$. This analysis is provided in Figure \ref{fig:Taylor-Green_surrogate_error}, depicting the maximum relative errors in $L^2(\mathcal{I}, H^1(\Omega))$, $L^\infty(\mathcal{I}, L^\infty(\Omega))$ and the time-gradient norm versus the reduced basis size. It shows a nearly exponential error decay as $N_\varepsilon$ increases. The maxima are computed on an independent test set, $\Xi_\text{TEST}^\mu \coloneqq \{ 1/(9.75 + 0.5 s)$, for all $ s \in \mathbb{N} \cap [1,80] \}$. Furthermore, the reduced solutions do not appear to deviate significantly from the projection of their full order counterparts onto the associated RB space, and the error bound employed demonstrates a good effectivity. 

For the implementation of the EnKM as presented in Subsection~\ref{sec:enkm}, it is necessary to provide a mathematical model for the measurement process. 
We take 40 measurements in time at the three sensor locations, $\eta_i$, $i \in \{1,2,3\}$, shown in Figure \ref{fig:domain_Taylor-Green}.
For this purpose, we introduce the measurement operator $\mathcal{L} : L^2(\mathcal{I}, L^2(\Omega)) \rightarrow \myR^{120}$, which can be seen as a vector of linear functionals $\ell_k : L^2(\mathcal{I}, L^2(\Omega)) \rightarrow \myR$ for all $k \in \myN \cap [1, 120]$. Each of those linear functionals has a unique Riesz representer $\rho_k : \mathcal{I} \times \Omega \rightarrow \myR$, with respect to the $L^2(\mathcal{I}, L^2(\Omega))$ norm, that can be written as
$$\rho_k = \nu_j \cdot \eta_i\quad\mbox{with}\quad k=3j+i\quad\mbox{for all}\: j\in \mathbb{N}\cap [1,40], \, i \in \mathbb{N}\cap [1,3],$$
where the spatial fields $\eta_i : \myGeoDomain \rightarrow \myR$ are Wendland functions $\psi_{2,1}$ of radius $0.1$ and center coordinates $(x_i,y_i) \in \{(0.1,0.7), (-0.1,-0.5), (0.5,0.1)\}$ (see Figure \ref{fig:domain_Taylor-Green}), while, for each $j\in\mathbb{N} \cap [1, 40]$, $\nu_j : \mathcal{I} \rightarrow \myR$ is a piecewise linear function supported over the interval $\mathcal{I}_j \coloneqq [t_j-2\Delta t, t_j+2\Delta t]$, where $t_j \coloneqq \Delta t (33+5j)$; $\nu_j$ is assumed to be symmetric with respect to $t_j$ and constant between $t_j-\Delta t$ and $t_j+\Delta t$.

Given this description of the observation process and the surrogate model, we next test the data assimilation scheme. We start with the estimation of the unknown parameter $\mu^\star = 0.04$ given the experimental measurements $\mathbf{y}(\mu^\star, \boldsymbol{\eta}) \in \myR^{120}$, with noise $\boldsymbol{\eta} \sim \mathcal{N}(\mathbf{0},\,\Sigma)$. We compare the performances of the EnKM employing a full order model and a surrogate model of accuracy $\varepsilon_c = 10^{-3}$ with $N_{\varepsilon}=42$. In order to obtain reliable statistics, we consider 25 ensembles $\mathcal{E}_0$ of size $J=150$ with particles sampled from the uniform prior distribution, $\Pi_0(\mu) = U(0.02, 0.10)$. The results obtained for a fixed value of $\sigma^2 = 10^{-6}$, at different iterations of the algorithm, are shown in Table \ref{tab:Taylor-Green_results}. We observe a quick stabilization of the error means $H_h$, $H_\varepsilon$ and $H_\varepsilon^*$, and of the error covariances, $S_h$, $S_\varepsilon$ and $S_\varepsilon^*$, after just a few steps. The full order algorithm performs significantly better than the biased reduced basis algorithm, while the adjusted version of the algorithm exhibits an excellent performance, very close to the full order one.

\begin{table}[t]
\caption{Comparison of reference FE $(\,\cdot_h)$ - biased RB $(\,\cdot_\varepsilon)$ - adjusted RB $(\,\cdot_\varepsilon^*)$ EnKM in low-noise conditions $\sigma^2=10^{-6}$. The test was performed by averaging 25 ensembles of 150 particles each and using reduced basis models of size $N_\varepsilon = 42$ ($\varepsilon_c \approx 0.001$). $H$ refers to the mean of the estimation error, while $S$ denotes the standard deviation of the estimation error, and c.t. the computational time.} 
\begin{center}
\begin{tabular}{c||c c || c c || c c}
        Iter & $H_h$ & $S_h$ & $H_\varepsilon$ & $S_\varepsilon$ & $H_\varepsilon^*$ & $S_\varepsilon^*$ \\[4pt]
        \hline 
        $0$ & $1.962$e-2 & $2.168$e-3 & $1.927$e-2 & $1.740$e-3 & $2.021$e-2 & $1.931$e-3  \\
        $1$ & $1.549$e-6 & $1.039$e-6 & $2.510$e-4 & $2.584$e-6 & $1.655$e-5 & $1.737$e-6  \\
        $2$ & $1.035$e-7 & $6.405$e-8 & $1.641$e-5 & $2.336$e-7 & $8.685$e-7 & $1.335$e-7  \\
        $3$ & $7.273$e-8 & $5.152$e-8 & $1.600$e-5 & $1.317$e-7 & $8.300$e-7 & $6.766$e-8  \\
        $4$ & $5.827$e-8 & $3.835$e-8 & $1.585$e-5 & $1.048$e-7 & $8.379$e-7 & $6.469$e-8  \\
        $5$ & $4.301$e-8 & $2.840$e-8 & $1.578$e-5 & $8.631$e-8 & $8.249$e-7 & $5.630$e-8  \\
        \hline
        c.t. & \multicolumn{2}{c ||}{ $2\text{h}\;\,54'\;\,10''$} &  \multicolumn{2}{c ||}{$3'\;\,13''$} &  \multicolumn{2}{c}{$3'\;\,07''$} 
\end{tabular}
\end{center}
\label{tab:Taylor-Green_results}
\end{table}

\begin{figure}[!b]
    \begin{center}
        \includegraphics[scale=0.72]{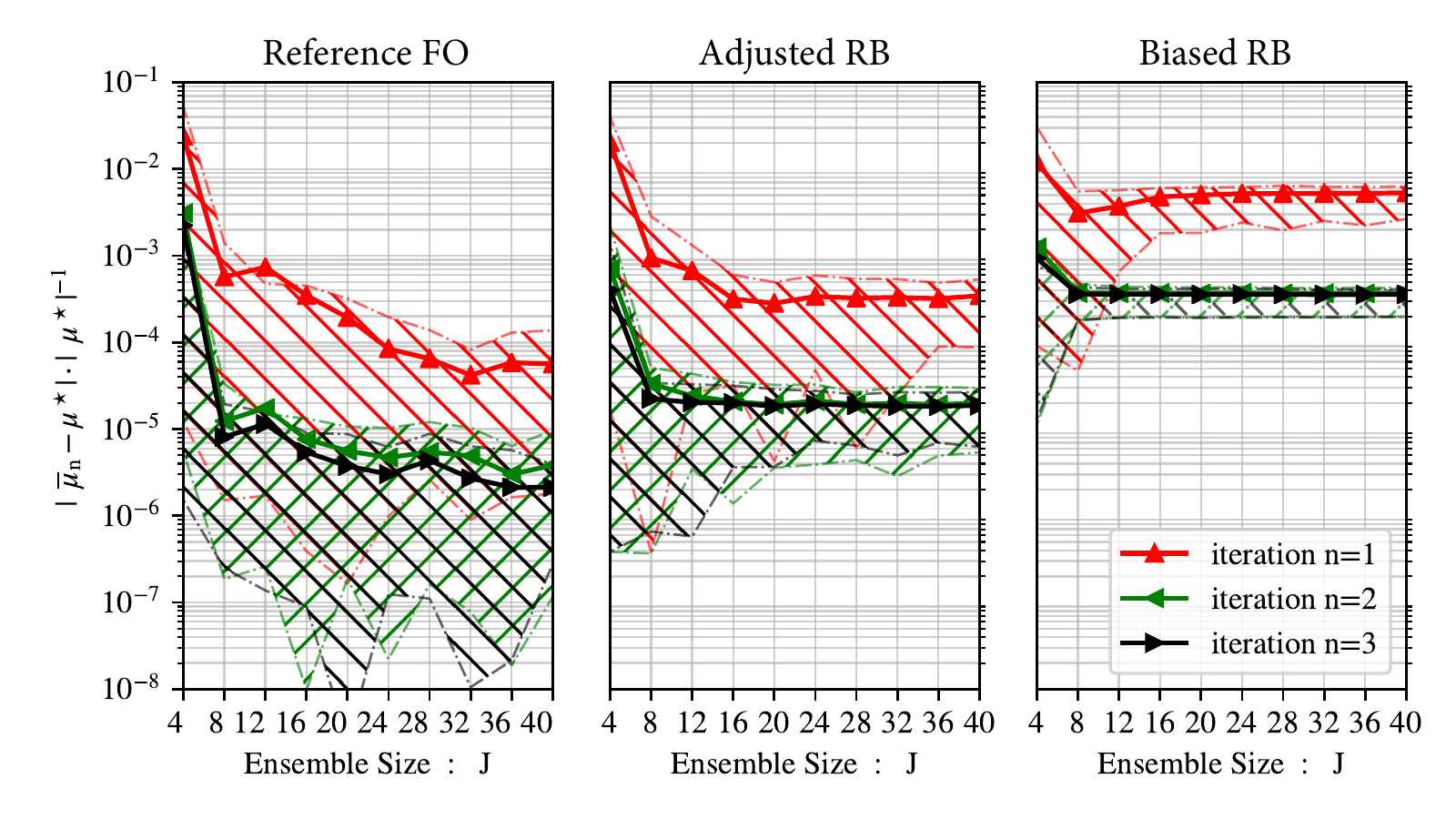}
    \end{center}
    \caption{Relative error in the parameter estimation vs. ensemble size $J$ for fixed noise magnitude, $\sigma = 10^{-3} \|\mathcal{G} (\myParaTrue)\|_\infty$. The standard full order EnKM is shown on the left, the adjusted RB-EnKM in the center and the biased RB-EnKM on the right. The solid lines represent the average error over 64 ensembles, while the dashed lines correspond to the 10th and 90th percentiles.}
    \label{fig:Taylor-Green_size_error}
\end{figure}

We next investigate the sensitivity of the algorithm to the accuracy of the reduced model, to the effect of the ensemble size, and to the noise magnitude. First, we repeat the estimation of the reference parameter $\mu^\star = 0.04$ for different values of the ensemble size $J = 4 k$, with $k \in \myN \cap [1, 10]$. In this experiment, we employ the same surrogate model used before and consider the relative noise magnitude $\sigma / \|\mathcal{G}(\myParaTrue)\|_\infty = 10^{-3}$. The results, shown in Figure \ref{fig:Taylor-Green_size_error}, indicate a larger sensitivity to $J$ for the full order algorithm than for the other two. It requires a larger number of particles before stabilizing on a large-ensemble asymptotic behavior (or mean-field behavior), while the reduced basis algorithms exhibit a much faster convergence, probably as a consequence of a lower-dimensional state space. Among the three iterations considered, the first appears to be the most affected, while, as the algorithm converges, the ensemble size seems to become less relevant.

\begin{figure}[!b]
    \begin{center}
        \includegraphics[scale=0.72]{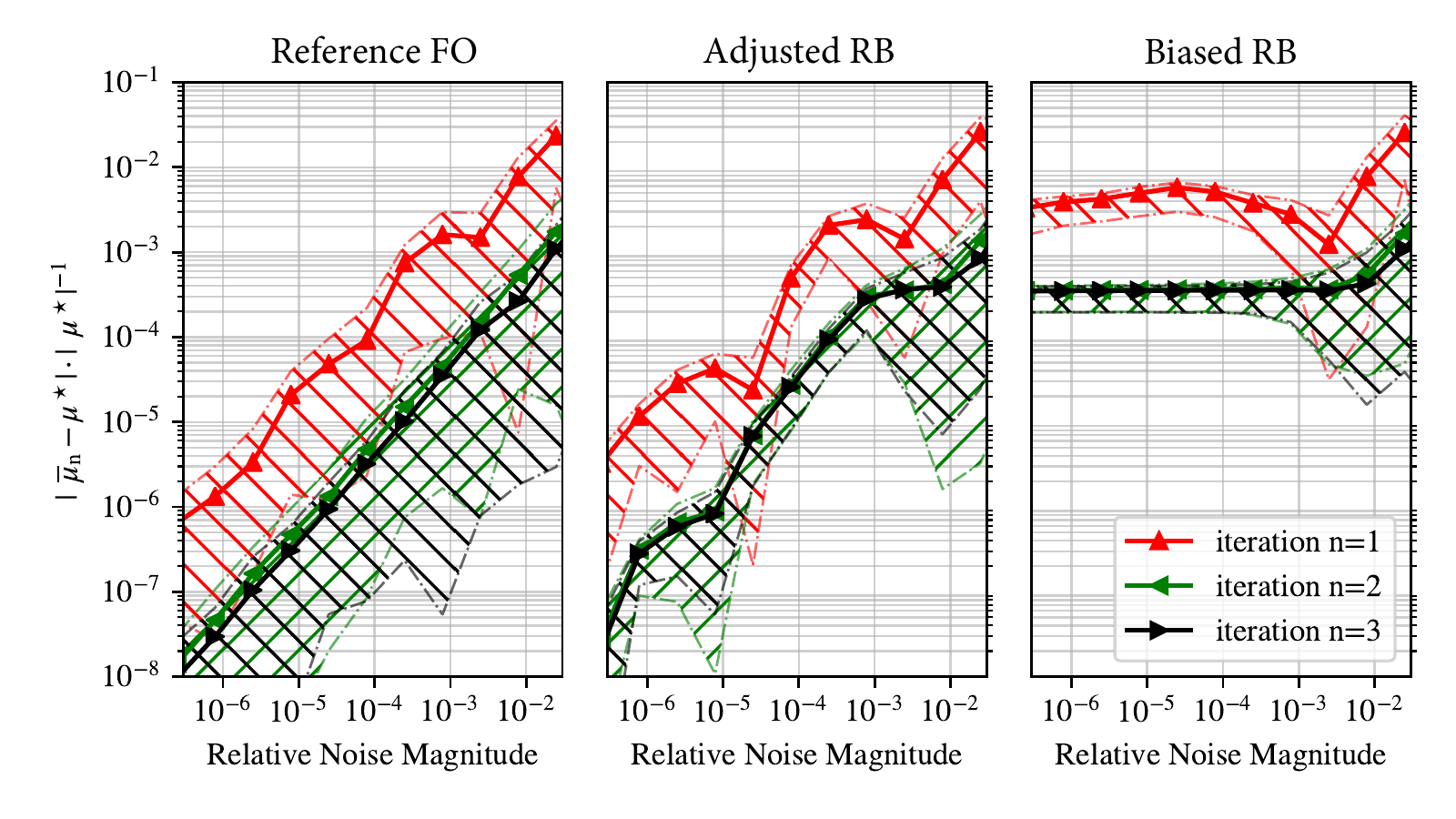}
    \end{center}
    \caption{Relative error in the parameter estimation vs. relative noise magnitude $\sigma /\|\mathcal{G} (\myParaTrue)\|_\infty$ for fixed ensemble size $J = 40$. The standard full order EnKM is shown on the left, the adjusted RB-EnKM in the center, and the biased RB-EnKM on the right. The solid lines represent the average error over 64 ensembles, while the dashed lines correspond to the 10th and 90th percentiles.}
    \label{fig:Taylor-Green_noise_error}
\end{figure}

In a second experiment, we consider the same parameter estimation, but we let the relative noise $\sigma /\|\mathcal{G}(\myParaTrue)\|_\infty$ take values $10^{-i}$ for $i\in\mathbb{N}\cap[2,6]$. Moreover, we employ $J=40$ particles per ensemble and the same reduced basis model as before. Each estimation is replicated 64 times for different noise realizations. The results are shown in Figure \ref{fig:Taylor-Green_noise_error}: for the full order EnKM we observe a linear dependence between the reconstruction error and the experimental noise, while the results for the biased RB-EnKM show that an untreated model bias introduces a systematic error independent of the noise magnitude. The most important result is the one related to the adjusted RB-EnKM: the error behavior achieved with this algorithm is comparable with the one obtained using a full order model. This demonstrates the effectiveness of the proposed method in compensating for the bias introduced by the reduced basis model, at least in this case study.

\begin{figure}[!t]
    \begin{center}
        \includegraphics[scale=0.75]{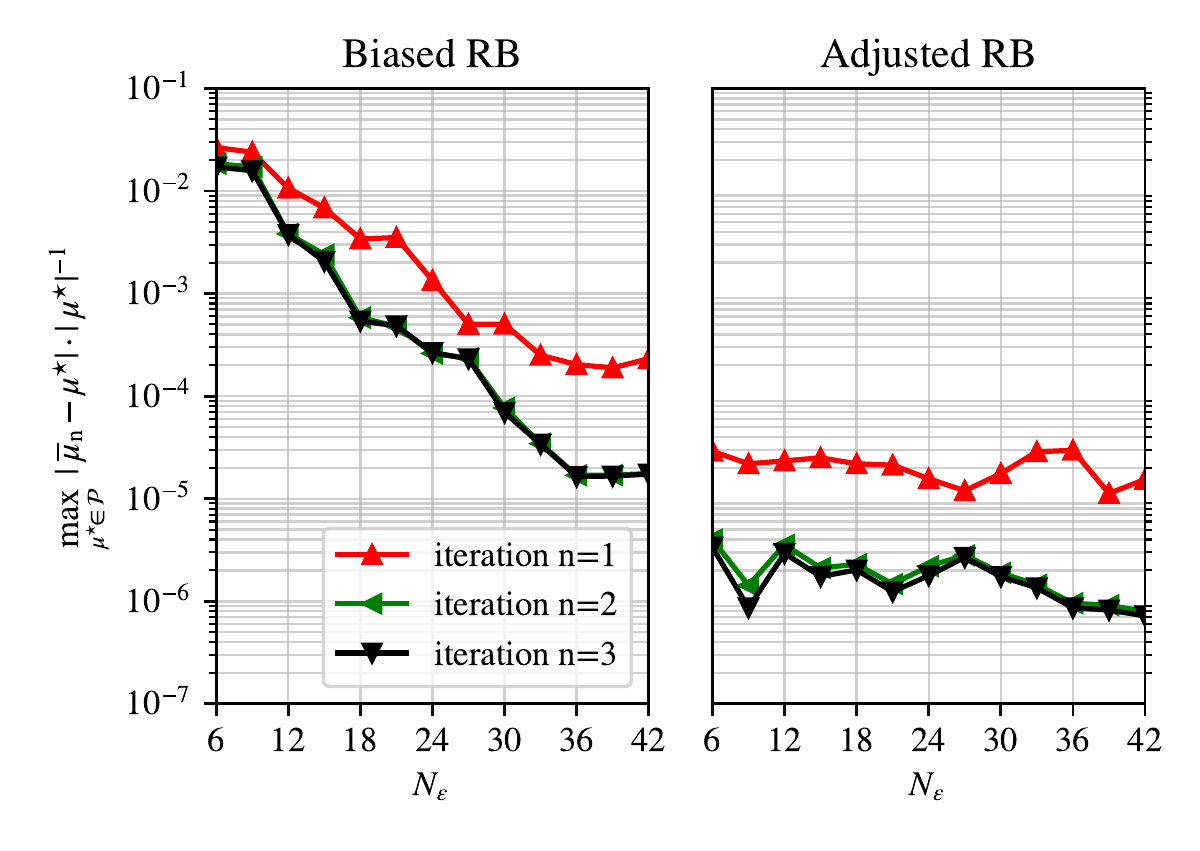}
    \end{center}
    \caption{Parameter error vs. reduced basis size for the biased and the adjusted RB-EnKM.}
    \label{fig:Taylor-Green_rb_relative_error}
\end{figure}

This conclusion is further confirmed by the last experiment, in which the performances of the biased and adjusted RB-EnKM are tested for all the parameters in the test set $\Xi_{\text{TEST}}^{\myPara}$ already employed to test the reduced basis model. Each parameter in the set is estimated using surrogate models of increasing sizes. Each estimation is performed 64 times in very low-noise conditions, that is $\sigma /\|\mathcal{G}(\myParaTrue)\|_\infty = 10^{-5}$, employing $J=40$ particles per ensemble. For each surrogate model employed, the results from the 64 ensembles are averaged and the maximum over the test set is computed. The results, shown in Figure \ref{fig:Taylor-Green_rb_relative_error}, demonstrate the ability of the correction to compensate for the presence of a model bias very well. As a consequence, the worst-case reconstruction error for the adjusted RB-EnMK barely depends on the reduced model size and it is always significantly lower than its biased counterpart. These results confirm the good performance of the adjusted RB-EnKM and its superiority over the biased RB-EnKM.

\subsection{Tracer Transport Problem}
\label{sec:nonlinear_experiment}

\begin{figure}[!t]
    \begin{center}
        \includegraphics[scale=0.7]{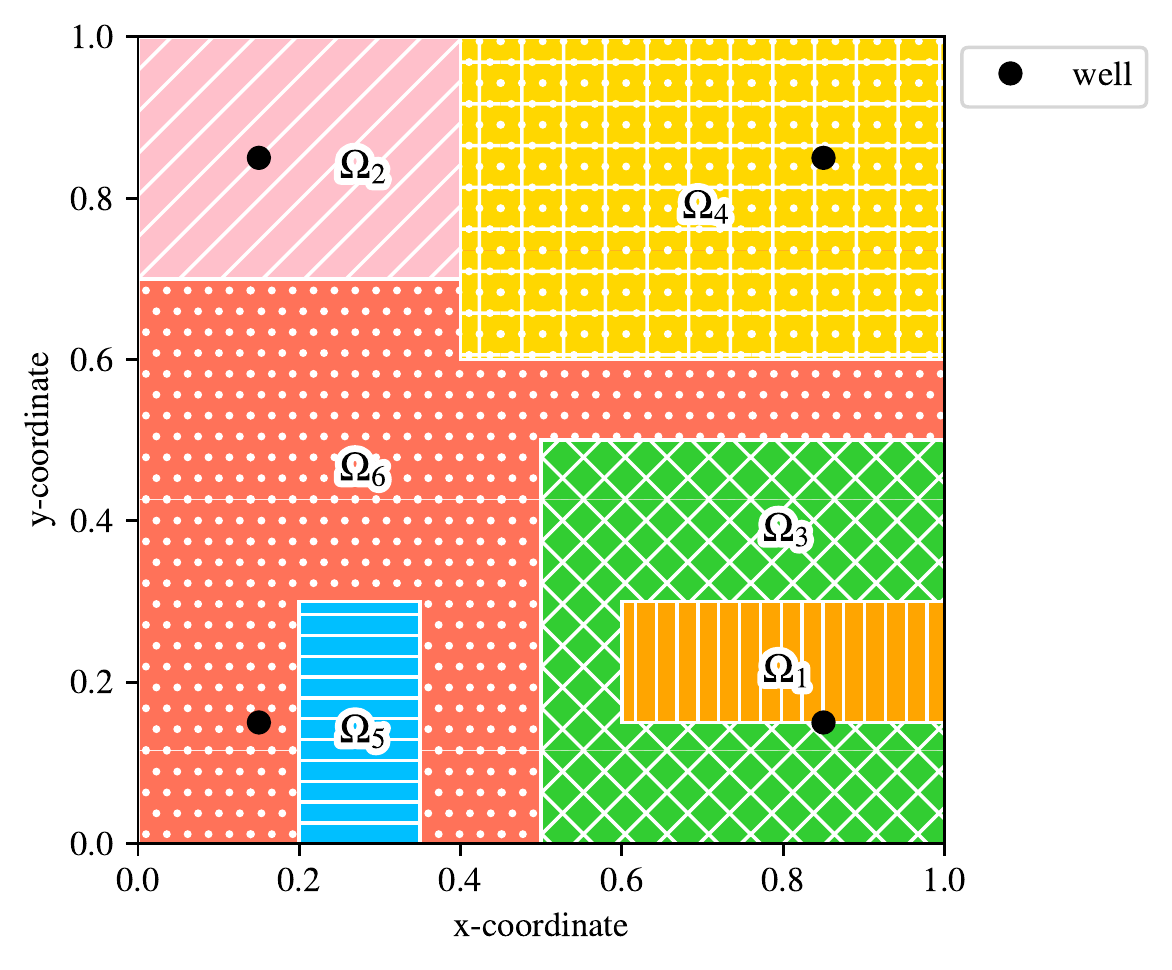}
        \caption{Domain of the tracer transport problem and injection wells.}
        \label{fig:domain_tracer}
    \end{center}
\end{figure}

We now consider the tracer transport problem from \cite{conrad2017parallel}, describing the non-homogeneous and non-isotropic transport of a non-reactive tracer in an unconfined aquifer. We introduce the spatial domain $\myGeoDomain \coloneqq (0, 1)^2$ divided into six sub-regions $\myGeoDomain = \bigcup_{r=1}^6 \myGeoDomain_r$ illustrated in Figure \ref{fig:domain_tracer} and defined as follows: $(x,y) \in \myGeoDomain$ is in $\myGeoDomain_r$ if the subscript $r$ is the smallest integer for which $x_0^r < x < x_1^r$ and $y_0^r < y < y_1^r$ where the points $\{(x_0^r,y_0^r)\}_{r=1}^6$ and $\{(x_1^r,y_1^r)\}_{r=1}^6$ are defined in Table \ref{Table1}.
We denote by $\partial \myGeoDomain$ the outer boundary of the domain and define the parallel walls $\Gamma_D \coloneqq (0,1)\times\{0,1\}$ and $\Gamma_N \coloneqq \partial \myGeoDomain \setminus \Gamma_D$. Based on this partition, we define the conductivity field as the piecewise constant function $k(\,\cdot\,;\myPara) : \myGeoDomain \rightarrow \myR$ over the six sub-regions $\myGeoDomain_r$. The conductivity can be affinely decomposed employing the coefficient vector $\myPara \in \myR^6$, with components $\mu_r$, and the indicator functions $\eta_r : \myGeoDomain \rightarrow \myR$
\begin{equation}
\label{eq:hydraulic_conductivity}
\begin{aligned}
    k(\myGeoCoord;\myPara) = \sum_{r=1}^6 e^{\mu_r} \eta_r(\myGeoCoord)  \qquad \text{with} \,\, \eta_r(\myGeoCoord) = \begin{cases} 
        1 \;\; \mbox{if} \quad \myGeoCoord \in \Omega_r, \\ 
        0 \;\; \mbox{if} \quad \myGeoCoord \in \Omega  \setminus \Omega_r. 
    \end{cases}
\end{aligned}
\end{equation}

We can now estimate the hydraulic log-conductivity $\myPara$, restricted to the orthotope $\mathcal{D} \coloneqq \varprod_{{}^{r=1}}^{{}_{6}} (\mu_r^{\min}, \mu_r^{\max}) \subset \myR^6$, relying on measurements of the tracer concentration $c(\,\cdot\,,\,\cdot\,;\myPara)$. This field satisfies the pPDE: find $c(\,\cdot\,,\,\cdot\,;\myPara) : \myGeoDomain \times (0,T] \rightarrow \myR$ such that
\begin{equation}
\label{eq:continuous_tracer_transport}
\left\{
\begin{alignedat}{3}
    &\partial_t c - \nabla \cdot ( ( d_m \mathbf{I} + d_l \boldsymbol{\beta} \boldsymbol{\beta}^\top ) \nabla c ) + \boldsymbol{\beta} \cdot \nabla c = f_c, \qquad
        && \mbox{in} \; \myGeoDomain \times \, (0,T]\,,\\
    & \nabla c(\mathbf{x},t;\myPara)\cdot \mathbf{n} = 0, 
        && \mbox{on} \;  \partial \myGeoDomain \times \, (0,T] \,,\\
    & c(\mathbf{x},0;\myPara) = 0, 
        && \mbox{in} \; \myGeoDomain\,.
\end{alignedat}
\right.    
\end{equation}
In this equation, the dispersion coefficients $d_l=d_m=2.5\cdot10^{-3}$ correspond to the flow-dependent component of the dispersion tensor and to its residual component, respectively. The forcing term $f_c$ is assumed to be of the form  $f_c \coloneqq \sum_{{}^{i=1}}^{{}_4} f_{c,i}$ and it models the injection of different amounts of tracer in four wells located at $(a_i, b_i)\in \{ 0.15, 0.85 \}^2$; each $f_{c,i}$ is a Gaussian function centered in $(a_i,b_i)$, with covariance $\Gamma_c = 0.005$ and multiplicative coefficient $p_i$ where $(p_1,p_2,p_3,p_4)=(10, 5, 10, 5)$. The velocity field $\boldsymbol{\beta} (\,\cdot\,,\,\myPara) : \myGeoDomain \rightarrow \myR^2$ is linearly dependent on the hydraulic head $u(\,\cdot\,,\,\myPara) : \myGeoDomain \rightarrow \myR$ through the relation $\boldsymbol{\beta} = -k(\myPara) \nabla u$. 
The latter field must satisfy the second constraint of the inverse problem, i.e., under the Dupuit–-Forchheimer approximation \cite{delleur2016} it solves the non-linear elliptic pPDE: find $ u(\,\cdot\,;\myPara) : \myGeoDomain \rightarrow \myR$ such that
\begin{equation}
\label{eq:continuous_hydraulic_head}
\left\{
\begin{alignedat}{3}
    & \nabla \cdot (k(\myPara) u \nabla u) + f_u = 0, \qquad
        && \mbox{in} \; \myGeoDomain \,,\\
    & \nabla u(\mathbf{x};\myPara) \cdot \mathbf{n} = 0, 
        && \mbox{on} \; \Gamma_N \,,\\
    &  u(\mathbf{x};\myPara)=0, 
        && \mbox{on} \; \Gamma_D \,.
\end{alignedat}
\right.    
\end{equation}
Here, the forcing term $f_u \coloneqq \sum_{i=1}^{4} f_{u,i}$ models an active pumping action at the four wells, each $f_{u,i}$ is a Gaussian function centered in $(a_i, b_i)$, of covariance $\Gamma_u = 0.02$ and coefficient $q_i$, where $(q_1,q_2,q_3,q_4)=(10, 50, 150, 50)$. Due to the combination of the quadratic dependence in $u$ and the zero boundary conditions, the equation always admits pairs of opposite solutions $u^{+}, u^{-}$. However, in our study, we are only interested in the positive solution $u^{+}(\;\cdot\,;\;\myPara) : \myGeoDomain \rightarrow \myRplus $.

\begin{figure}[!b]
    \begin{center}
        \includegraphics[scale=0.78]{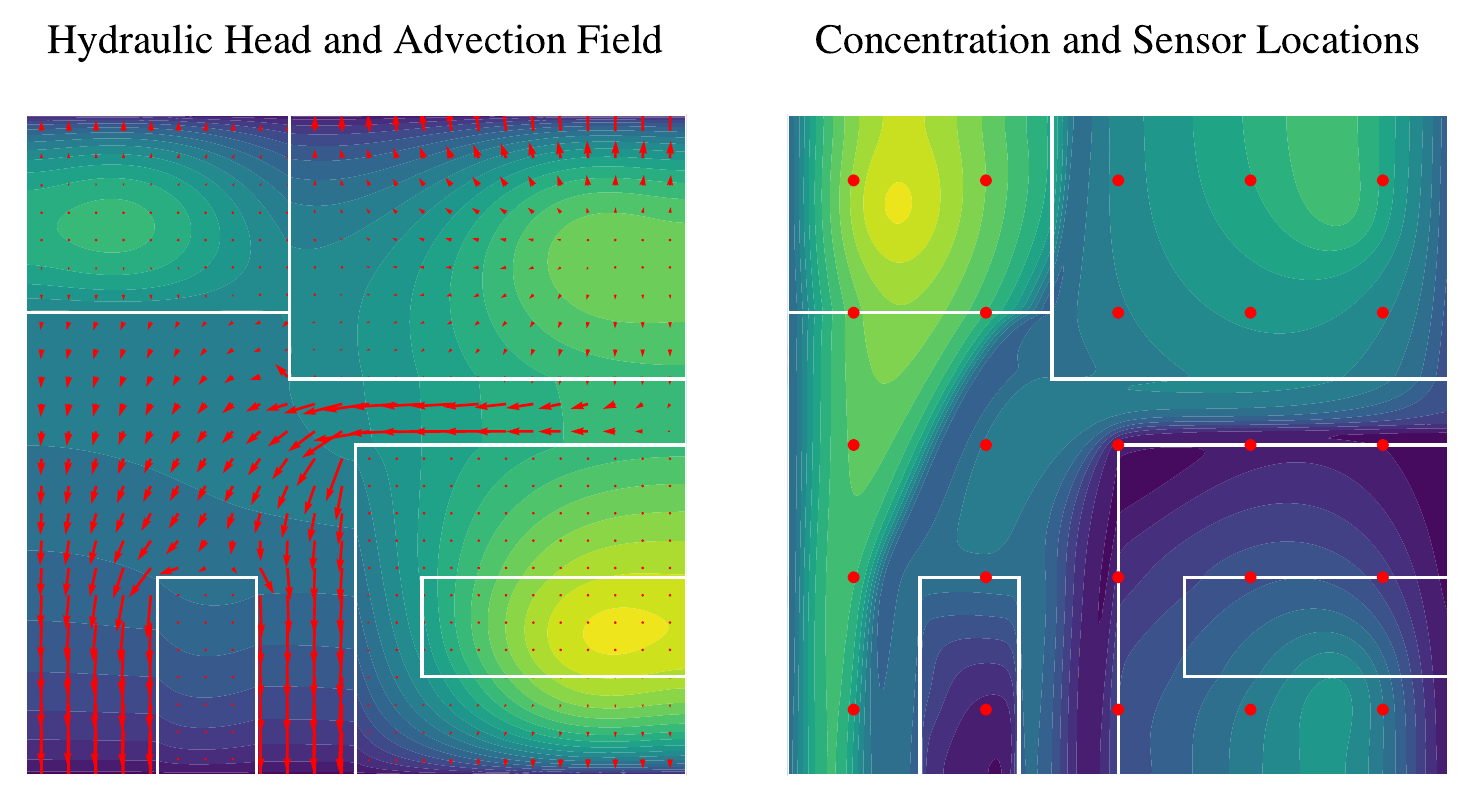}
        \caption{Reference solutions at log-conductivity ${\myParaTrue} := [-0.75, -0.25, -0.5, \;1, -0.25, \;3]$. On the left: hydraulic head $u_h(\myParaTrue)$ and corresponding velocity field $\boldsymbol{\beta}_h\coloneqq - k(\myParaTrue) \nabla u_h$ (in red). On the right: tracer concentration $c_h(\myParaTrue)$, at time $t=0.4$, and measurement wells (in red).}
        \label{fig:ref_solution_tracer}
    \end{center}
\end{figure}

\begin{table}[t]
\caption{On the left: coordinates of the corners of the sub-regions $\Omega_r$. On the right: true values of the parameters $\mu_r$ and boundaries of the uniform prior $\Pi_0$} 
\begin{center}
\begin{tabular}{c|c c c c ||c c c | c}
        Region & $\quad x_0^r$ & $x_1^r$ & $y_0^r$ & $y_1^r\quad$ & ($\star$) Ref & Min & Max $\quad$ & $\quad \mu_r$ \\[4pt]
        \hline 
        $\Omega_1$ & $\quad 0.60$ & $1.00$ & $0.15$ & $0.30\quad$ & $\;\; -0.75$ & $-1.00$ & $0.00\quad$ & $\quad \mu_1$ \\
        $\Omega_2$ & $\quad 0.00$ & $0.40$ & $0.70$ & $0.10\quad$ & $\;\; -0.25$ & $-1.00$ & $1.00\quad$ & $\quad \mu_2$ \\
        $\Omega_3$ & $\quad 0.50$ & $1.00$ & $0.00$ & $0.50\quad$ & $\;\; -0.50$ & $-1.00$ & $0.00\quad$ & $\quad \mu_3$ \\
        $\Omega_4$ & $\quad 0.40$ & $1.00$ & $0.60$ & $1.00\quad$ & $\;\; \;\;1.00$ & \;\;$0.00$ & $2.00\quad$ & $\quad \mu_4$ \\
        $\Omega_5$ & $\quad 0.20$ & $0.25$ & $0.00$ & $0.30\quad$ & $\;\; -0.25$ & $-1.00$ & $0.00\quad$ & $\quad \mu_5$ \\
        $\Omega_6$ & $\quad 0.00$ & $1.00$ & $0.00$ & $1.00\quad$ & $\;\; \;\;3.00$ & \;\;$2.00$ & $5.00\quad$ & $\quad \mu_6$ \\
        \hline
\end{tabular}
\end{center}
\label{Table1}
\end{table}

Full order solutions are obtained via a finite element approximation, employing piecewise linear functions, $\zeta_i: \myGeoDomain \rightarrow \myR$, for $i = 1, \ldots, N_h$, with $N_{h} = 44,972$ degrees of freedom (mesh size $h \approx 0.01$). The discretization of the elliptic equation \eqref{eq:continuous_hydraulic_head} results in a discrete non-linear problem which is iteratively solved employing a Newton scheme with tolerance $10^{-6}$. The approximate solution, $u_{h}$, is used to compute the velocity field, $\boldsymbol{\beta}_h\coloneqq - k(\myPara) \nabla u_{h}$, which is piecewise constant with $N_{h} - 1$ degrees of freedom. This is needed for the solution of the parabolic equation \eqref{eq:continuous_tracer_transport}, whose discretization leads to a system of ordinary differential equations integrated over the time interval $\mathcal{I} = \left(0.0,\,0.5\right]$ using the Crank--Nicolson scheme with uniform time step $\Delta t = 0.01$.
This is equivalent to performing a Petrov--Galerkin projection of Equation \eqref{eq:linear_advection_diffusion} with trial and test spaces defined as follows: we consider the partition of the temporal interval $\mathcal{I}$ into the union of equispaced subintervals $\mathcal{I}_n \coloneqq \left(t_{n-1}, t_n \right]$ of length $\Delta t$ with $n=1,\ldots, N_t$ and $N_t \coloneqq 0.5 / \Delta t$.
Let $\omega_n :\mathcal{I} \rightarrow \myR$ be a piecewise constant function with support in $\mathcal{I}_n$, and let $\upsilon_n :\mathcal{I} \rightarrow \myR$ be a hat function with support in $\mathcal{I}_{n} \cup \mathcal{I}_{n+1}$. We define the trial space $\mathcal{V}_h \coloneqq \text{span} \{ \upsilon_n \cdot \zeta_i \}_{ {}^{i,n=1} }^{ {}_{N_h, N_t} }$ and the test space $\mathcal{W}_h \coloneqq \text{span} \{ \omega_n \cdot \zeta_i \}_{ {}^{i,n=1} }^{ {}_{N_h, N_t} }$, respectively .

Each full order simulation is obtained employing a FreeFEM++ solver \cite{hecht2012new} and takes roughly 2 minutes to be computed. Figure \ref{fig:ref_solution_tracer} shows the hydraulic head $u_h(\,\cdot\,;\myParaTrue)$ and the relative velocity field $\boldsymbol{\beta}_h(\,\cdot\,;\myParaTrue)$ (on the left) and the tracer concentration field $c_h(\,\cdot\,,0.4;\myParaTrue)$ (on the right), both associated with the reference log-conductivity 
\begin{align}
    \myParaTrue = [-0.75, -0.25, -0.50, 1.00, -0.25, 3.00]^\top. 
    \label{eq:ref_para}
\end{align}
The same reference log-conductivity is used as the true parameter for the data assimilation problem. Pointwise observations are collected at five successive times $t_m \in \{0.1, 0.2, 0.3, 0.4, 0.5\}$, in $25$ spatial location $\mathbf{x}_{ij} = (x_i, y_j)$ such that $x_i=0.1+0.2 i$ and $y_j=0.1+0.2 j$ for $i,j \in \{0, \ldots,4\}$. This operation is encoded in the measurement operator $\mathcal{L} : H^1(\myGeoDomain) \rightarrow \myR^{125}$. Each noise-free measurement is polluted with i.i.d. Gaussian noise with mean zero and covariance $\sigma$, resulting in a noise covariance matrix $\Sigma = \sigma^2 \mathbf{I}$. 

In order to solve the inverse problem with surrogate models of different accuracy, many approximations of equation \eqref{eq:continuous_hydraulic_head} and \eqref{eq:continuous_tracer_transport} must be produced. This requires the introduction of spatial basis functions $\psi_i, \varphi_j : \myGeoDomain \rightarrow \myR$, $i \in \myN \cap [1, N_\varepsilon]$, $j \in \myN \cap [1, M_\varepsilon]$, selected by applying the method of snapshots (POD) to the two sets of full order solutions, 
\[\Theta^u_\text{TRAIN} \coloneqq \{ u_h (\,\cdot\,;\myPara^{{}_{(s)}}) \}_{{}^{s=1}}^{{}_{S}}\quad\text{and}\quad\Theta_\text{TRAIN}^{c} \coloneqq \{ c_h(\,\cdot\,, t^{\scriptscriptstyle  (z)};\myPara^{{}_{(s)}}) \}_{{}^{z,s=1}}^{{}_{Z,S}},\]
with snapshot parameters, $\myPara^{{}_{(s)}} \sim \Pi_0 \coloneqq \varprod_{r=1}^6 U(\mu^{\text{min}}_r, \mu^{\text{max}}_r)$, for all $s \in \myN \cap [1, S]$ and sampling times $t^{\scriptscriptstyle (z)} = 0.01 z$ for all $z \in \myN \cap [1, Z]$, where $S=2,000$ and $Z=50$. The number of basis functions considered, $N_\varepsilon, M_\varepsilon \in \myN$, is the one required to approximate the hydraulic head and the tracer concentration with relative accuracy $\varepsilon_u, \varepsilon_c \in \myRplus$, where 
\begin{align}
    \label{eq:hydraulic_head_accuracy_definition}
    \varepsilon_u &\coloneqq \sup_{\myPara \in \mathcal{D}} \frac{\| u_h(\myPara) - u_\varepsilon(\myPara) \|_{H^1(\myGeoDomain)}}{\|u_h(\myPara)\|_{H^1(\myGeoDomain)}}, \\ 
    \label{eq:tracer_concentration_accuracy_definition}
    \varepsilon_c &\coloneqq \sup_{\myPara \in \mathcal{D}} \frac{\| c_h(\myPara) - c_\varepsilon(\myPara) \|_{L^2(\mathcal{I}, H^1(\myGeoDomain))}}{\|c_h(\myPara)\|_{L^2(\mathcal{I}, H^1(\myGeoDomain))}}.
\end{align}

\begin{figure}[!b]
    \begin{center}
        \includegraphics[scale=0.72]{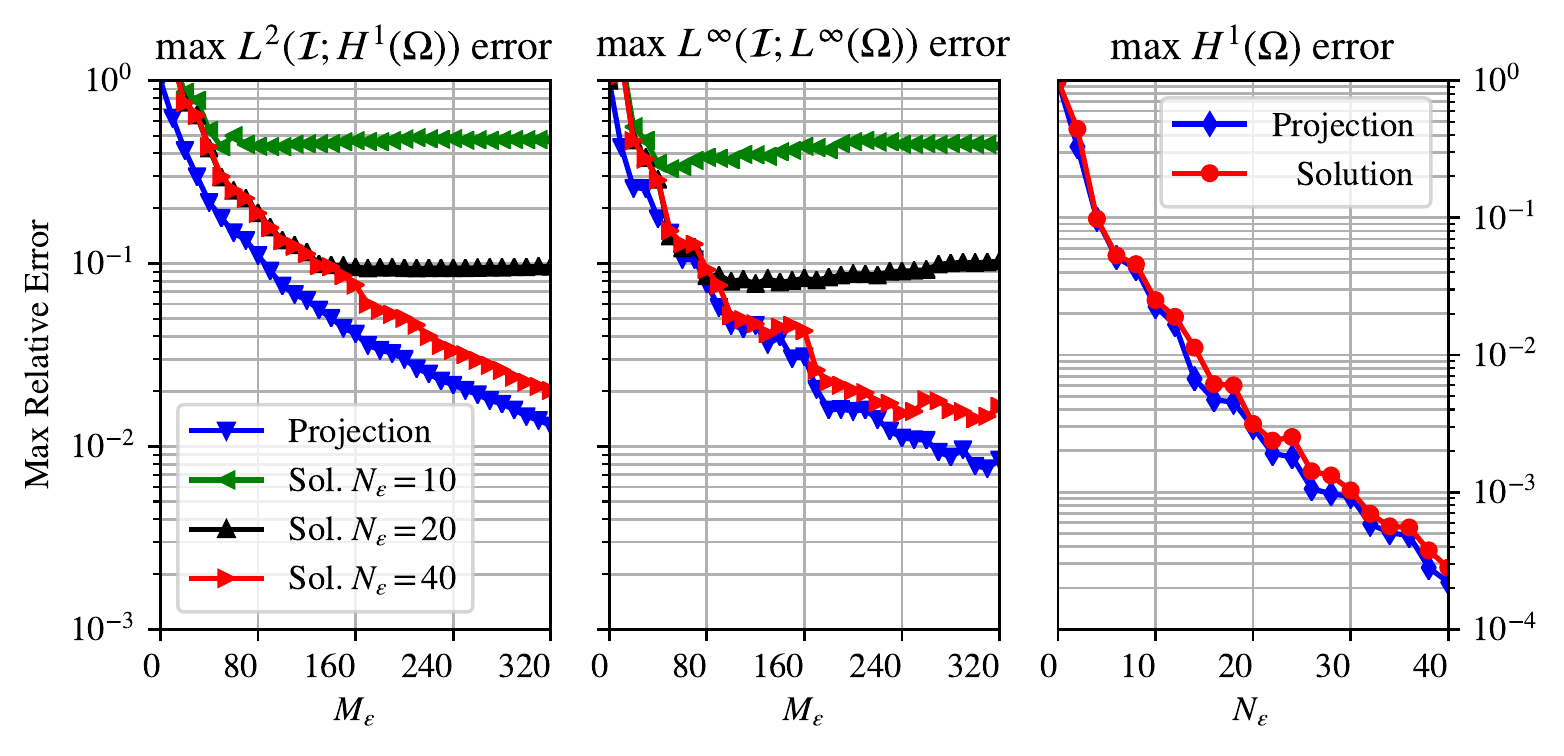}
        \caption{Left and center: maximum relative error of the solution and of the projection of the tracer concentration vs. $M_\varepsilon$ for different values of $N_\varepsilon$; projection -- in space -- performed with respect to the $H^1(\Omega)$ inner product. Right: maximum relative error of the projection and of the solution of the hydraulic head vs. $N_\varepsilon$. Error norm shown above each plot.}
        \label{fig:tracer_transport_surrogate_error}
    \end{center}
\end{figure}

Based on the first set of basis functions, the approximation space for the Galerkin projection of equation \eqref{eq:continuous_hydraulic_head} is defined as $\mathcal{U}_\varepsilon \coloneqq \text{span} \{ \psi_i \}_{{}^{i=1}}^{{}_{N_\varepsilon}}$. From the second set of basis functions, instead, the RB test space $\mathcal{W}_\varepsilon \coloneqq \text{span} \{ \omega_n \cdot \varphi_i \}_{ {}^{i,n=1} }^{ {}_{M_\varepsilon, N_t} }$ and RB trail space $\mathcal{V}_\varepsilon \coloneqq \text{span} \{ \upsilon_n \cdot \varphi_i \}_{ {}^{i,n=1} }^{ {}_{M_\varepsilon, N_t} }$ are defined for the Petrov--Galerkin projection of equation \eqref{eq:continuous_tracer_transport}.
We look at reduced solutions of the form 
\begin{align}
    u_\varepsilon (\mathbf{x};\myPara) &= \sum_{i=1}^{N_\varepsilon} u_{i}(\myPara)\psi_i(\myGeoCoord)&\text{ for } \myGeoCoord \in \myGeoDomain,\\
    c_\varepsilon (\mathbf{x}, t;\myPara) &= \sum_{j=1}^{M_\varepsilon} \sum_{n=1}^{N_t} c_{n,j}(\myPara) \upsilon_n(t) \varphi_j(\myGeoCoord) &\text{ for } t\in\mathcal{I},\,\myGeoCoord \in \myGeoDomain.
\end{align}
Where the expansion coefficients  $c_{n,j}$ and $u_i$, with $i\in \myN \cap [1, N_\varepsilon]$, $n \in \myN \cap [1,N_t]$ and $j \in \myN \cap [1,M_\varepsilon]$, satisfy the systems of algebraic equations %
\begin{align}
    \label{eq:alg_hydraulic_head}
    \sum_{p,q=1}^{N_\varepsilon, N_\varepsilon} \mathbf{N}_{ipq} (\myPara) u_p u_q  &= f_i, \\ 
    \label{eq:alg_tracer_concentration}
    \sum_{k=1}^{ M_\varepsilon} \left( \mathbf{M}_{jk} + \frac{\Delta t}{2} \mathbf{D}_{jk} (\mathbf{u},\myPara) \right) c_{n+1,k} &= \left(  \mathbf{M}_{jk} - \frac{\Delta t}{2} \mathbf{D}_{jk} (\mathbf{u},\myPara) \right) c_{n,k} + g_j,
\end{align}
given the initial conditions $c_{0,j}=0$ for all $j \in \myN \cap [1, M_\varepsilon]$. The scalar forcing terms $f_i$, $g_j$ are obtained by integrating their full order counterparts versus the basis functions $\psi_i$ and $\varphi_j$, for all $i \in \myN \cap [1, N_\varepsilon]$, $j \in \myN \cap [1, M_\varepsilon]$
\begin{align}
    f_i \coloneqq \int_\Omega f_h \psi_i d \Omega, \qquad g_j \coloneqq \Delta t \int_\Omega f_c \varphi_j d \Omega.
\end{align}
The mass and stiffness matrices $\mathbf{M}, \mathbf{K} \in \myR^{\scriptscriptstyle  M_\varepsilon}$ are defined as in \eqref{eq:matrices}, while the parameter dependent tensors $\mathbf{D}(\mathbf{u},\myPara) \in \myR^{\scriptscriptstyle  M_\varepsilon^2}$ and  $\mathbf{N}(\myPara) \in \myR^{\scriptscriptstyle  M_\varepsilon^3}$ depend affinely on the multidimensional arrays $\mathbf{A}\in \myR^{\scriptscriptstyle 6 \times N_\varepsilon^3}$, $\mathbf{B}\in \myR^{\scriptscriptstyle  6 \times N_\varepsilon^2 \times M_\varepsilon^2}$, and $\mathbf{C}\in \myR^{\scriptscriptstyle  6 \times N_\varepsilon \times M_\varepsilon^2}$ defined as
\begin{align}
    \label{eq:rb_tensor_A}
    \mathbf{A}_{ipqr} &\coloneqq \int_\Omega \frac{\eta_r}{2} \left( \psi_p (\nabla \psi_q \cdot \nabla \psi_i) + \psi_q (\nabla \psi_p \cdot \nabla \psi_i) \right) d\Omega, \\
    \label{eq:rb_tensor_B}
    \mathbf{B}_{jkpqr} &\coloneqq \int_\Omega \eta_r (\nabla \varphi_j \cdot \nabla \psi_p)(\nabla \varphi_k \cdot \nabla \psi_q) d\Omega, \\
    \label{eq:rb_tensor_C}
    \mathbf{C}_{jksr} &\coloneqq \int_\Omega \eta_r (\nabla \varphi_j \cdot \nabla \psi_s) \varphi_k d\Omega. 
\end{align}
For a fixed value of the log-conductivity, $\myPara$, the tensors $\mathbf{N} (\myPara) $ and $\mathbf{D} (\mathbf{u}, \myPara) $ can be assembled. The latter, however, requires the evaluation of the discrete hydraulic head $\mathbf{u}$. They are respectively defined as
\begin{align}
    \label{eq:rb_tensor_N}
    \mathbf{N}_{ipq} (\myPara) &\coloneqq \sum_{r=1}^{6}  e^{\mu_r} \mathbf{A}_{ipqr} \,, \\
    \label{eq:rb_tensor_D}
    \mathbf{D}_{jk} (\mathbf{u}, \myPara) &\coloneqq d_m \mathbf{K}_{jk} + d_l \sum_{p, q, r=1}^{N_\varepsilon, N_\varepsilon, 6} e^{2\mu_r} \mathbf{B}_{jkpqr} u_p u_q  + \sum_{s,r=1}^{N_\varepsilon, 6} e^{\mu_r} \mathbf{C}_{jksr} u_s \,.
\end{align}
We emphasize that the accuracy of the solutions to Equations \eqref{eq:alg_hydraulic_head} and \eqref{eq:alg_tracer_concentration}, with the latter equivalent to a Crank--Nicolson discretization, depends on the number of basis functions and on the time step $\Delta t$. In Figure \ref{fig:tracer_transport_surrogate_error}, on the left and on the right, we show the maximum relative errors of the surrogate model (${\varepsilon}_u$, ${\varepsilon}_c$) as a function of $N_\varepsilon$ and $M_\varepsilon$. In the center, we show the $L^\infty( \mathcal{I} ;  L^\infty( \myGeoDomain ))$ relative error of the tracer concentration, bounding from above the error on synthetic measurements. We compute these maximum relative errors on a set of parameters $ \Xi_\text{TEST}^{\myPara} \coloneqq \{ \myPara^{{}_{(s)}} \sim \Pi_0(\myPara) \}_{s=1}^{500}$ independent of the ones used for the model training. It can be observed that, for small values of $N_{\varepsilon}$, the error in the concentration stagnates after a certain value of $M_{\varepsilon}$, suggesting that, in this region, the error is dominated by the approximation of the hydraulic head. However, for $N_\varepsilon=40$, this effect is no longer present, at least for the values of $M_\varepsilon$ considered, and the tracer error only depends on $M_{\varepsilon}$. This allows us to modify the accuracy of the model by varying the dimension of the reduced model.

The construction of the reduced model has an offline cost of about 100 hours. This includes the time required for the construction of a training set of $2,000$ full order solutions $(61\text{h}\, 16'\, 40'')$, the time for the computation of the POD basis functions $(23' \, 40'')$, and the time for assembling the RB model tensors $(13\text{h}\, 32'\, 47'')$. This cost corresponds roughly to the computational cost of $2,500$ finite element solutions, each of which takes approximately $110$s. By contrary, the surrogate model obtained employing $N_\varepsilon = 40$, $M_\varepsilon = 320$ basis functions produces a solution in only $1.25$s (online cost), which is about $1/90$ of its full order equivalent. The same training set used for the POD is employed to estimate, at negligible cost, the empirical moments of $\boldsymbol{\delta}^\star_\varepsilon$, i.e., $\overline{\boldsymbol{\delta}}_\varepsilon$ and $\boldsymbol{\Gamma}_\varepsilon$. 

We now focus on the inverse problem, as discussed in Section \ref{sec:enkm}. We start by considering the estimation of the reference parameter $\myParaTrue$ given the measurements $\mathbf{y}(\myParaTrue, \boldsymbol{\eta}) \in \myR^{125}$, polluted by experimental noise of magnitude $\sigma$. To have a reliable statistic, we consider $32$ independent initial ensembles $\mathcal{E}_{0}$ of variable size, sampled from the same distribution $\Pi_0$. 

As a first experiment, we compare the performances of the two RB-EnKM employing $J=160$ particles and a surrogate model with error tolerance $\varepsilon_c \approx 0.02$ (obtained with $N_\varepsilon = 40$ and $M_\varepsilon = 320$). The first test relies on the biased version of the RB-EnKM, as presented in Section \ref{sec:enkm}, while the second test corresponds to the adjusted algorithm. For both simulations, we consider low-amplitude experimental noise, i.e., negligible if compared to the model error, $\sup_{\myPara} \| \mathcal{L} (c_h(\myPara) - c_\varepsilon(\myPara)) \|_\infty \approx 10^{-2} > 10^{-3} = \sigma$, and we separately pollute the measurements employed by the ensembles. In Table \ref{tab:EnKF_ref_RB_vs_RB-adj}, we report the average properties of the ensembles after $4$ iterations: columns $E_\varepsilon$ and $E_\varepsilon^*$ contain the mean parameter estimation, i.e., the particle mean, averaged over the $32$ ensembles. Here, columns $\Sigma_\varepsilon$ and $\Sigma_\varepsilon^*$ contain the average standard deviation of the ensembles. We can observe that the correction term has the effect of significantly lowering the reconstruction error from $\| E_\varepsilon - \myParaTrue \|_\infty = 6.437$e-$3$ to $\| E_\varepsilon^* - \myParaTrue \|_\infty = 7.870$e-$4$.  We also notice that the variability of the estimate increases consistently with the presence of an additional term in the Kalman gain.
\begin{table}[t]
\caption{Comparison of biased RB $(\,\cdot_\varepsilon)$ - adjusted RB $(\,\cdot_\varepsilon^*)$ EnKM in low-noise conditions $\sigma = 10^{-6}$. The test was performed by averaging $32$ ensembles, each employing $160$ particles and $4$ iterations, and using reduced basis models of size $N_\varepsilon=40$, $M_\varepsilon=320$ ($\varepsilon_c \approx 0.02$). $E$ refers to the average parameter estimation, while $\Sigma$ denotes the average ensemble standard deviation, and c.t. the computational time.} 
\begin{center}
\begin{tabular}{c || c c || c c || c}
        $          $ & $\;\; E_{\varepsilon}$ & $\Sigma_{\varepsilon}$ & $\;\;E_{\varepsilon}^*$ & $\Sigma_{\varepsilon}^*$ & $\quad \myParaTrue$ \\[4pt]
        \hline  
        $\mu_1 \;\;$ & $  - 0.751890$ & $1.261\text{e-3} \;\;$ & $  - 0.750787$ & $1.543\text{e-3} \;\;$ & $\;  - 0.75\;$ \\
        $\mu_2 \;\;$ & $  - 0.256437$ & $1.860\text{e-3} \;\;$ & $  - 0.250403$ & $2.785\text{e-3} \;\;$ & $\;  - 0.25\;$ \\
        $\mu_3 \;\;$ & $  - 0.503134$ & $1.001\text{e-3} \;\;$ & $  - 0.500397$ & $1.353\text{e-3} \;\;$ & $\;  - 0.50\;$ \\
        $\mu_4 \;\;$ & $\;\;0.999042$ & $0.701\text{e-3} \;\;$ & $\;\;0.999764$ & $1.016\text{e-3} \;\;$ & $\;\;\;1.00\;$ \\
        $\mu_5 \;\;$ & $  - 0.251110$ & $0.898\text{e-3} \;\;$ & $  - 0.250325$ & $1.221\text{e-3} \;\;$ & $\;  - 0.25\;$ \\
        $\mu_6 \;\;$ & $\;\;3.000127$ & $0.630\text{e-3} \;\;$ & $\;\;2.999665$ & $0.904\text{e-3} \;\;$ & $\;\;\;3.00\;$ \\
        \hline
        $\text{c.t.}$ \;\; & \multicolumn{2}{c ||}{ $7\text{h}\;\,12'\;\,59''$} & \multicolumn{2}{c ||}{ $7\text{h}\;\,11'\;\,42''$} & $ \;1'\;\,51''$\\
\end{tabular}
\end{center}
\label{tab:EnKF_ref_RB_vs_RB-adj}
\end{table}

\begin{figure}[!t]
    \begin{center}
        \includegraphics[scale=0.78]{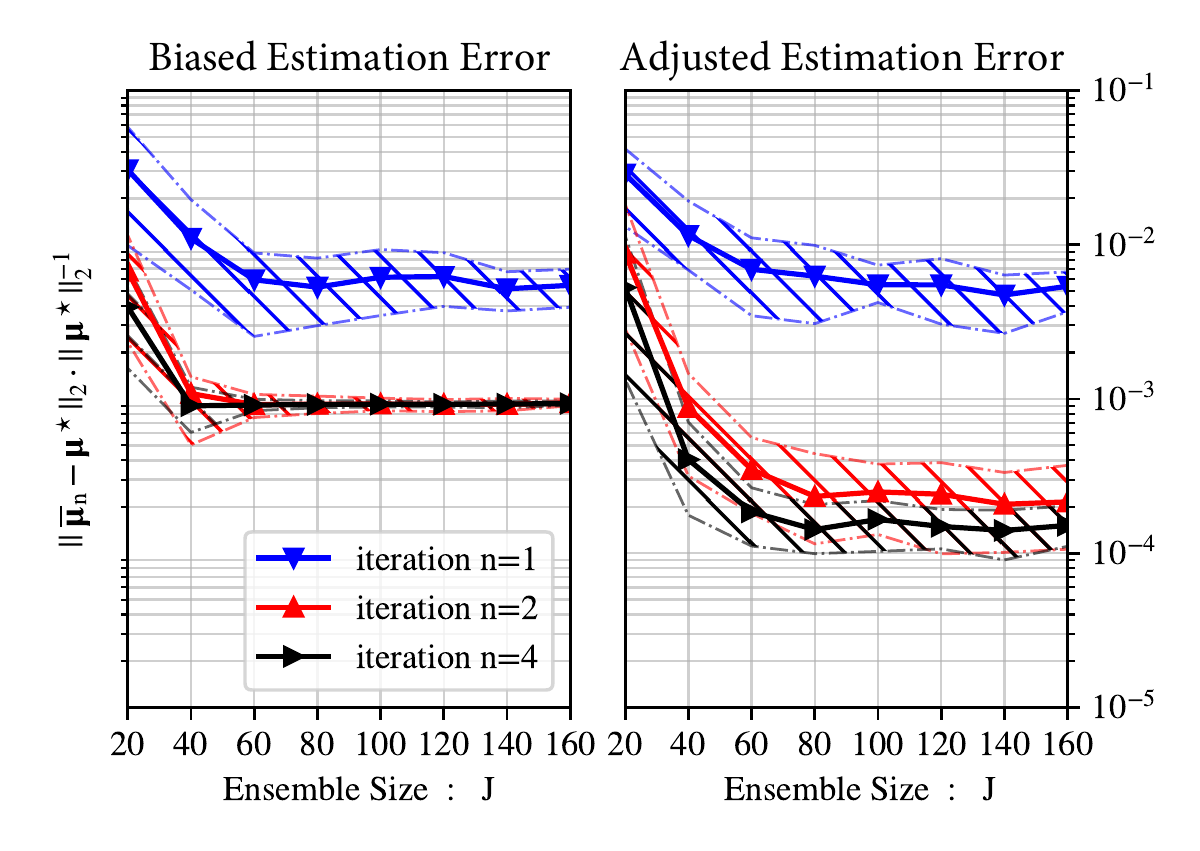}
        \caption{Biased and adjusted RB-EnKM parameter estimation relative error vs. ensemble size, for fixed noise magnitude $\sigma = 10^{-3}$. The solid lines represent the average error over 32 ensembles at different algorithm iterations. The dashed lines represent the 10th and the 90th percentiles.}
        \label{fig:tracer_en_size_error_comparison}
    \end{center}
\end{figure}

As an extension of the previous experiment, we estimate the reference parameter $\myParaTrue$ employing the same surrogate model, noise magnitude and number of ensembles as before, but using ensembles of variable size $J=20k$, with $k \in \myN \cap [2,16]$. This allows us to study the effect of the ensemble size on the parameter estimation obtained with the biased and adjusted RB-EnKM algorithms. The results shown in Figure \ref{fig:tracer_en_size_error_comparison} indicate that, for both algorithms, very small ensembles lead to large relative errors and entail a large variability among the different samples. This behavior seems to be relevant only for ensembles with less than $40$ particles when the biased RB-EnKM is employed, and with less than $80$ particles when the adjusted version is used. Larger ensembles do not exhibit relevant fluctuations; we can therefore assume an ensemble of size $J=160$ to be sufficiently large to ensure the independence from this quantity of the results in the upcoming tests.

\begin{figure}[!t]
    \begin{center}
        \includegraphics[scale=0.78]{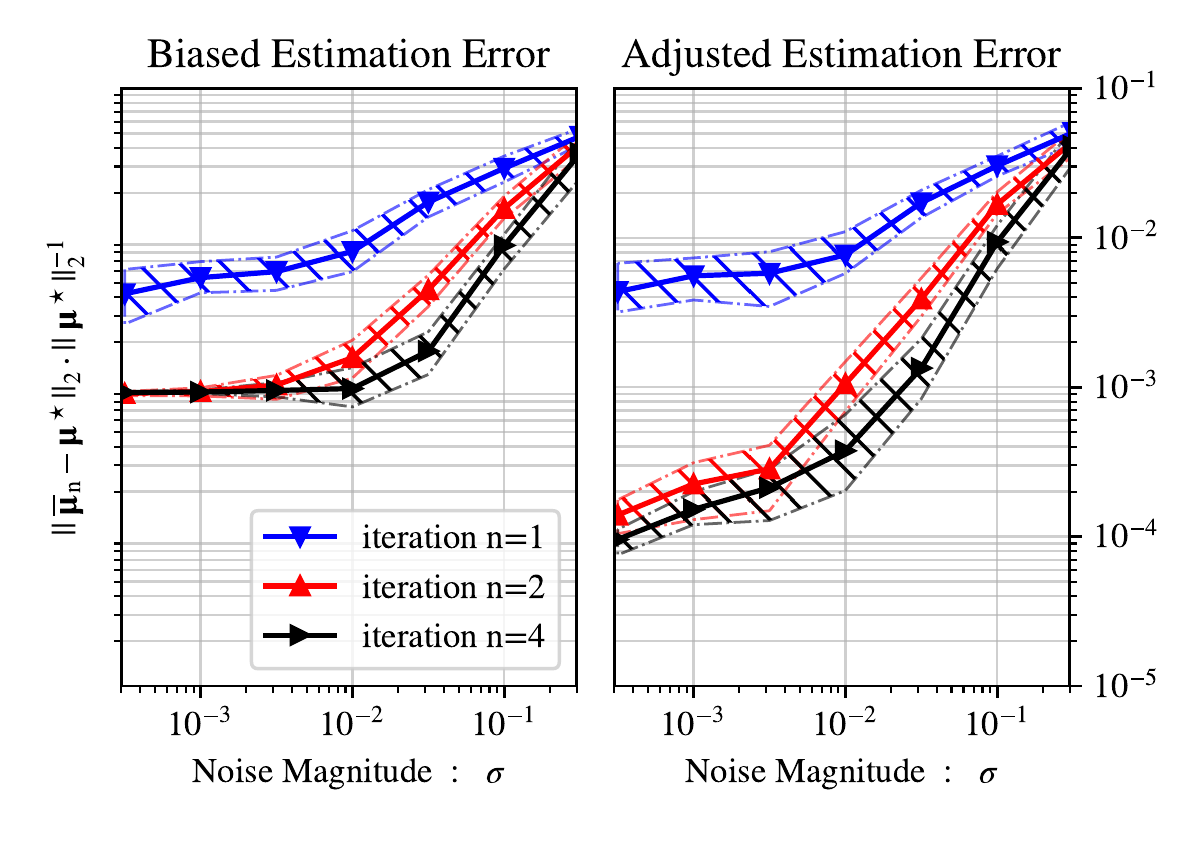}
        \caption{Biased and adjusted RB-EnKM parameter estimation relative error vs. absolute noise magnitude, for fixed ensemble size $J=160$. The solid lines represent the average error over 32 ensembles at different algorithm iterations. The dashed lines represent the 10th and the 90th percentiles.}
        \label{fig:tracer_noise_comparison}
    \end{center}
\end{figure}

A key quantity determining the performances of the method is the noise magnitude. Its effect on the two reduced basis algorithms is investigated by looking at the variation of the relative estimation error of the reference parameter $\myParaTrue$ when the noise magnitude varies. To this end, we consider seven noise values, $\sigma^2 = 10^{-m}$ with $m \in \myN \cap [1, 7]$. We employ the same RB-EnKM used before, with a fixed ensemble size $J=160$, and we average the results over $32$ independent ensembles. The results, shown in Figure \ref{fig:tracer_noise_comparison}, reiterate the inadequacy of the biased method in dealing with the systematic bias introduced in the measurements by the surrogate model. In fact, the plot corresponding to the biased method shows error stagnation for low-noise. On the contrary, the plot corresponding to the adjusted method highlights a mitigation of this effect, with an estimation error that keeps decreasing in low-noise conditions, although at a lower rate than in high-noise conditions.

\begin{figure}[!b]
    \begin{center}
        \includegraphics[scale=0.78]{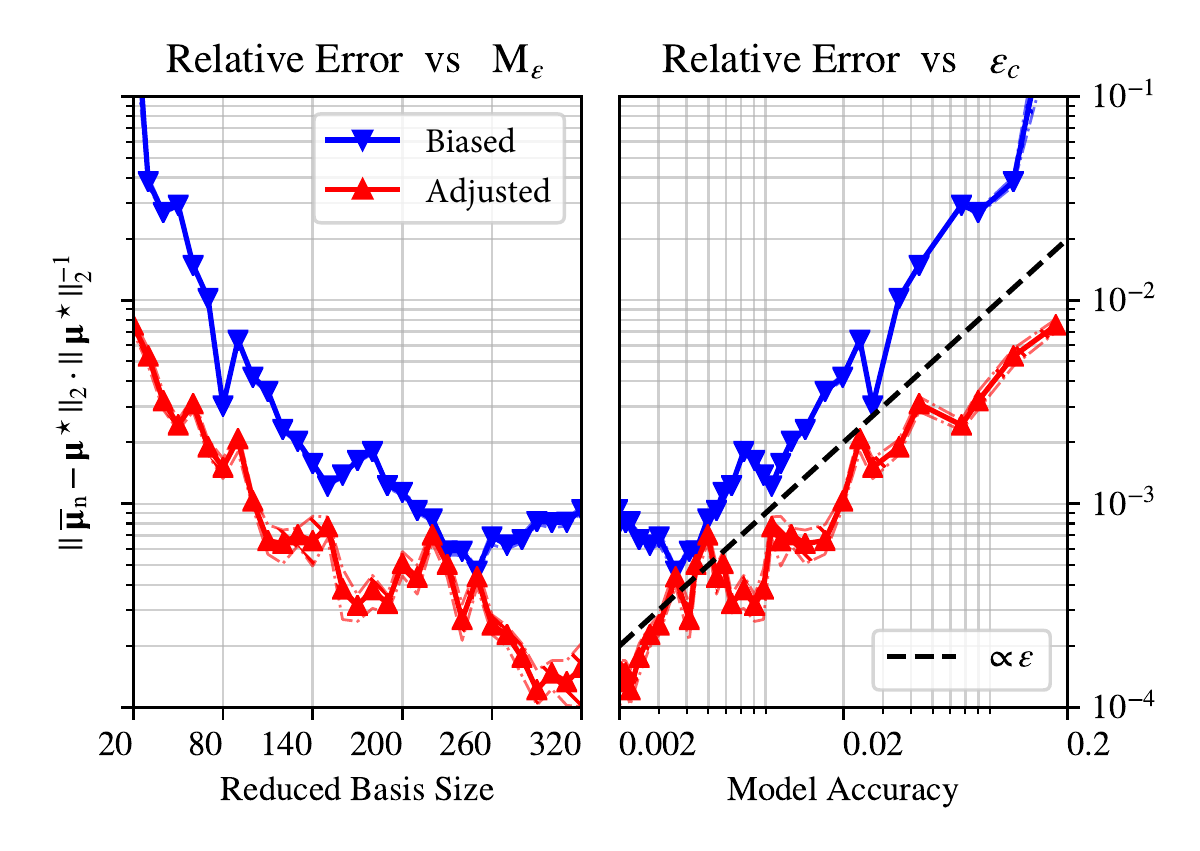}
        \caption{Parameter error vs. reduced basis size and maximum relative error of the solution. The solid lines represent the average error over 16 ensembles at different algorithm iterations. The dashed lines represent the 10th and the 90th percentiles.}
        \label{fig:tracer_rb_comparison}
    \end{center}
\end{figure}

In our last experiment, we test the performances of the biased and adjusted RB-EnKM by employing surrogate models of increasing accuracy. We fix the size of the reduced space $\mathcal{U}_\varepsilon$ to a sufficiently large value, $N_\varepsilon = 40$, and we vary the size of the approximation space associated with the concentration: $M_\varepsilon =10k$, with $k \in \myN \cap [2, 32]$. Employing the resulting approximated models, we estimate the reference parameter $\myParaTrue$ in low-noise conditions, $\sigma = 10^{-3}$, averaging the results obtained over $16$ ensembles of $160$ particles each. In Figure \ref{fig:tracer_rb_comparison}, we show the final relative error (after three algorithm iterations) as $M_\varepsilon$ and $\varepsilon_c$ change, both for the biased and the adjusted RB-EnKM. For both, we observe that the relative estimation error decreases, almost linearly, with the error of the surrogate model. Moreover, we observe that, with few exceptions, the error of the adjusted algorithm is smaller than the error of the biased algorithm. The few points where the two errors are very close can be explained by a strongly unbalanced distribution of the measurement bias in a region away from the reference parameter. Future developments that take into account, in the execution of the algorithm, the parameter estimate to adjust the bias correction should dampen this effect.

\section{Conclusions}
\label{sec:conclusions}

We proposed an efficient, gradient-free iterative solution method for inverse problems that combines model order reduction techniques, via the reduced basis method, and the Kalman ensemble method introduced in \cite{Iglesias_2013}. The use of surrogate models allows a significant speed-up of the computational cost, but it leads to a distortion in the cost function optimized by the inverse problem.  This in turn introduces a systematic error in the approximate solution of the inverse problem. To overcome this limitation, we have proposed the adjusted RB-EnKM which corrects for this bias by systematically adjusting the cost function and thus retrieving good convergence.

Using a linear Taylor--Green vortex problem, the performance of the method is compared versus the full order model as well as to the biased RB-EnKM in which no adjustment was made. The numerical results show that the biased method fails to achieve the same accuracy of the full order method. Contrarily, the adjusted RB-EnKM attains the same accuracy as its full order counterpart for a large range of noise magnitudes at a significantly lower computational cost, and even approaches the mean-field limit faster as the ensemble size is increased. Furthermore, the dependence on model accuracy of the reconstruction error is essentially removed over the range of model accuracy considered. 

The method was then applied to a non-linear tracer transport problem for which the full order model was impractical. The results for this example show that, despite a decrease in the order of convergence at low-noise, the stagnation of the reconstruction error observed in the biased RB-EnKM can be removed by adjusting the algorithm. Regarding the model accuracy, a substantial improvement of the adjusted EnKM with respect to the biased EnKM was observed, although less pronounced than in the linear problem. 

Overall, our numerical tests show that the proposed method allows for the use of inexpensive surrogate models while empirically ensuring that the predicted result of the inversion remains accurate with respect to the full order inversion at a significantly lower computational cost.


\newpage
\bibliography{sn-bibliography}


\end{document}